\documentclass[11pt,a4paper,reqno]{amsart}

\usepackage{amsfonts,amsthm,amscd,epsfig,amsmath,amssymb,enumerate}

\numberwithin{equation}{section}

\DeclareMathSymbol{\leqslant}{\mathalpha}{AMSa}{"36} % nicer `smaller or equal'
\DeclareMathSymbol{\geqslant}{\mathalpha}{AMSa}{"3E} % nicer `larger or equal'
\DeclareMathSymbol{\eset}{\mathalpha}{AMSb}{"3F}     % nicer `emptyset'
\renewcommand{\leq}{\;\leqslant\;}                   % redef. of < or =
\renewcommand{\geq}{\;\geqslant\;}                   % redef. of > or =
\newcommand{\be}{\begin{equation}}

\def\1{\ifmmode {1\hskip -3pt \rm{I}} \else {\hbox {$1\hskip -3pt \rm{I}$}}\fi}

\newtheorem{Th}{Theorem}[section]
\newtheorem{Le}[Th]{Lemma}
\newtheorem{Pro}[Th]{Proposition}

\newtheorem{Def}[Th]{Definition}
\title[From combinatorics to large deviations]{From combinatorics to large deviations
for the invariant measures of some multiclass particle systems}

\author{Davide Gabrielli}
\address{Dipartimento di Matematica Universit\`a dell'Aquila, Via Vetoio Loc. Coppito 67100 L'Aquila,
Italy. E--mail: gabriell@univaq.it}

\begin{document}

\begin{abstract}
We prove large deviation principles (LDP) for the invariant measures
of the multiclass totally asymmetric simple exclusion process
(TASEP) and the multiclass Hammersely-Aldous-Diaconis (HAD) process
on a torus. The proof is based on a combinatorial representation of
the measures in terms of a \emph{collapsing procedure} introduced in
\cite{A} for the $2$-class TASEP and then generalized in \cite{FM1},
\cite{FM2} and \cite{FM3} to the multiclass TASEP and the multiclass
HAD process. The rate functionals are written in terms of
variational problems that we solve in the cases of $2$-class
processes.
\end{abstract}

\maketitle

\section{Introduction}

In recent years several new and interesting results have been
obtained in the study of fluctuations of interacting particle
systems. Some of these results concern the computations of large
deviations rate functionals for specific models.

Given a stochastic interacting particle system, a problem of
interest is the determination of its invariant measures. When the
model is not reversible and the detailed balance does not hold, this
can be a difficult task. Typical examples are boundary driven
stochastic lattice gases. Depending on the model you can have
available some representations of the invariant measures or not. We
will give in this introduction a short outline of some of the recent
progress in computation of the rate functionals of large deviations
for the empirical measures of the invariant measures in models of
this type. The results are interesting for several reasons. One
reason is that the measures have often long range correlations and
the corresponding rate functionals are not local. They have a
structure very different from the one obtained in the case for
example of Gibbs measures where you have an integration of a
function of the density profile. Another reason is that interacting
particle systems are very effective models of statistical mechanics
and the results obtained give insight for the behavior of more
complex models or real systems.

From one side there are combinatorial representations for the
invariant measures of exclusion like models starting from which it
is possible to compute the corresponding rate functionals (see
\cite{D} for a recent review and references therein). These
combinatorial representations are built up from products of
operators satisfying appropriate commutation relations.

A different approach is based on dynamical arguments. Fluctuations
of the invariant measure can be recovered from fluctuations of paths
of the processes. The static rate functional for the invariant
measure is then obtained solving a variational problem for the
dynamic rate functional. This leads to an Hamilton-Jacobi equation
as a central object (see \cite{BDGJL} for a recent review and
references therein). Differently from the exact solution approach
the dynamic one is insensitive to small perturbations of the
dynamics corresponding to the same macroscopic structure of
fluctuations. Nevertheless the Hamilton-Jacobi equation is, in
general, difficult to solve.

In this paper we prove large deviation principles (LDP) for the
invariant measures of the multiclass totally asymmetric simple
exclusion process (TASEP) and the multiclass
Hammersely-Aldous-Diaconis (HAD) process. Our proof is based on an
exact combinatorial representation of the invariant measures.
Configurations of particles distributed according to the invariant
measures are constructed applying a deterministic transformation,
the \emph{collapsing procedure}, to configurations of particles
distributed according to product of uniform measures. We generalize
the collapsing procedure up to let it act on positive measures. An
application of the contraction principle allows one to obtain the
final result. The rate functionals that we obtain are not local. In
the case of $2$ class models they are obtained from a geometric
construction on density profiles.

As our results are obtained from the contraction principle, the rate
functionals are naturally expressed as infimum of auxiliary
functionals. Several non local functionals obtained as rate
functionals of particle systems are represented in terms of either
infimum or supremum of auxiliary functionals. An interesting
question is whether there is always a representation of them as an
infimum. More precisely if those rate functionals can be obtained
from wider LDP using the contraction principles. This is the case of
the present paper. This is also the case of the invariant measures
for the TASEP with boundary sources as represented in \cite{DuSch}.
The corresponding application of the contraction principle has been
done in \cite{BDGJL}. In \cite{BDGJL} it is also suggested that this
could be the case for the KMP model with boundary sources. A rate
functional obtained from contraction of a convex rate functional is
not necessarily convex, this is also the case of the present result.

The paper is organized as follows. In section \ref{multimodel} we
define the multiclass TASEP and the multiclass HAD process
constructing them using the basic coupling. In section \ref{colla}
we describe the \emph{collapsing procedure} as it acts on
configurations of particles. We discuss also briefly its main
properties. In section \ref{invcolla} we show how the collapsing
procedure is used to construct the invariant measures of the
processes. In section \ref{collames} we generalize the collapsing
procedure defining its action on positive measures. We discuss also
its main properties. In section \ref{collaempmis} we define
empirical measures and relate them to the collapsing procedure. In
section \ref{ldpu} we derive, from well known results, large
deviation principles for uniform distributions. In section
\ref{ldp2} we derive LDP for the invariant measures of the 2-class
TASEP and the 2-class HAD process. Using the contraction principle
we have the rate functionals in a variational form for which we can
find the unique minimum in terms of a \emph{concave envelope}
construction. These are theorems \ref{Th1} and \ref{Th2} that are
the main results of the paper. The rate functionals are not convex.
We discuss also in variational terms the typical density of first
class particles when the total density is known and the typical
total density when the density of first class particles in known. In
section \ref{ldpk} we derive LDP for the multiclass TASEP and the
multiclass HAD process. The rate functionals are written in terms of
variational problems. It is interesting to study such problems and
others proposed in the section. We obtain also a recursive relation.

The combinatorial constructions of \cite{A}, \cite{FM1}, \cite{FM2},
\cite{FM3} are different from the original solution of the $2$-class
TASEP obtained in \cite{DJLS}, that is based on products of non
commuting operators. It is interesting to derive the same result of
the present paper starting from this alternative solution. A
probably different representation of the same rate functional will
maybe appear. It is also interesting to study the problem using the
dynamic approach. The major problem here is the lack of a complete
dynamical LDP. A problem of interest is also the study of the
variational problem \eqref{ultimo} for models with more than $2$
classes of particles.

To avoid confusion we remark that we use similar symbols for
mathematical objects that play a similar role in the TASEP, the HAD
process and in the general framework of positive measures. We use
also the same symbol $\mathbb C$ to indicate the \emph{collapsing
operator} both when it acts on configurations of particles and on
positive measures. This is due to the fact that the second one is a
natural generalization of the previous one.

\section{Multiclass models}
\label{multimodel}

\subsection{TASEP}

The totally asymmetric simple exclusion process (TASEP) is a model
of stochastic jumping particles satisfying an exclusion rule. Let
$\mathbb{Z}_N:=\mathbb{Z}/(N\mathbb Z)$ be the discrete one
dimensional torus with $N$ sites. Every site $x\in \mathbb{Z}_N$ can
be either empty or occupied by a particle. The state space of the
process is $X_N:=\left\{0,1\right\}^{\mathbb Z_N}$. Given $\eta\in
X_N$ a configuration of particles, we will say that the site $x\in
\mathbb Z_N$ is occupied by a particle if $\eta(x)=1$ and empty
otherwise. Every particle waits an exponential time of rate one and
then tries to jump to its nearest neighbor site to the left. If the
site is already occupied by another particle then the jump is
suppressed. This informal description can be summarized by the
following generator of the dynamics
\begin{equation}
L_Nf(\eta):=\sum_{x\in \mathbb
Z_N}\left(f(\eta^{e_x})-f(\eta)\right)
%\label{gentasep}
\nonumber
\end{equation}
where $e_x$ is the oriented bond $(x,x+1)$ and $\eta^{e_x}$ is the
configuration of particles obtained from $\eta$ rearranging the
values at the extremes of $e_x$ in decreasing order, according to
the following definition
\begin{equation}
\eta^{e_x}(z):=\left\{
\begin{array}{lcl}
\eta(z) & if & z\neq x,x+1\ , \\
\max \left\{\eta(x), \eta(x+1)\right\} & if & z=x\ , \\
\min \left\{\eta(x), \eta(x+1)\right\} & if & z=x+1 \ .\\
\end{array}
\right. \nonumber
\end{equation}

A generalization of the previous model is obtained labeling some of
the particles as \emph{first class particles} and the remaining ones
as \emph{second class particles}. When a first class particle tries
to jump over a second class particle it succeeds and the two
particles exchange their positions. When a second class particle
tries to jump over a first class particle the jump is suppressed.
The natural state space for such a process is
$\left\{0,1,2\right\}^{\mathbb Z_N}$ obtained from the choice of
assigning value $0$ to empty sites, value $1$ to sites occupied by a
first class particle and value $2$ to sites occupied by second class
particles. We will instead describe a configuration with a pair
$(\eta_1,\eta_2)$ with both $\eta_1$ and $\eta_2$ elements of $X_N$.
The configuration $\eta_1\in X_N$ is such that $\eta_1(x)=1$ when in
$x$ there is a first class particle and $\eta_1(x)=0$ otherwise. The
configuration $\eta_2\in X_N$ is such that $\eta_2(x)=1$ if in $x$
there is either a first or a second class particle and $\eta_2(x)=0$
otherwise. We endow $X_N$ with the natural partial order $\preceq$
defined from
$$
\eta\preceq \xi \ \ \Leftrightarrow \eta(x)\leq \xi(x), \ \ \forall
x\in \mathbb Z_N \ .
$$
By definition $\eta_1\preceq \eta_2$ so that $(\eta_1,\eta_2)\in
I_N^{2,\uparrow}$ where
$$
I_N^{k,\uparrow}:=\left\{(\eta_1,\dots ,\eta_k): \eta_i\in X_N, \
\eta_i\preceq \eta_{i+1}\right\}\ .
$$
The two descriptions of the state space are equivalent and a
bijection between $I_N^{2,\uparrow}$ and
$\left\{0,1,2\right\}^{\mathbb Z_N}$ is defined from
$$
\xi(x):=\inf\left\{i: \eta_i(x)=1, \ i=0,1,2\right\}\ ,
$$
where we defined $\eta_0(x):=1-\eta_2(x)$.

\noindent The generator of the above described 2-class TASEP is
$$
L_Nf(\eta_1,\eta_2):=\sum_{x\in \mathbb
Z_N}\left(f(\eta_1^{e_x},\eta_2^{e_x})-f(\eta_1,\eta_2)\right)\ .
$$
This generator clearly defines also a jointly, order preserving,
evolution of two TASEP, usually called basic coupling.

A further natural generalization, called the k-class TASEP is
obtained introducing particles of class up to a fixed natural number
$k\leq N$. When a particle of class $i$ tries to jump over a
particles of class $j$ with $j>i$ then the positions of the two
particles are exchanged. When a particle of class $i$ tries to jump
over a particle of class $j$ with $j\leq i$ the jump is suppressed.
The state space is now $I_N^{k,\uparrow}$ that is in bijection with
$\left\{0,1,\dots ,k\right\}^{\mathbb Z_N}$. The generator of the
dynamics is
$$
L_Nf(\eta_1,\dots ,\eta_k):=\sum_{x\in \mathbb
Z_N}\left(f(\eta_1^{e_x},\dots ,\eta_k^{e_x})-f(\eta_1,\dots
,\eta_k)\right)\ .
$$
This generator clearly defines also a jointly, order preserving,
evolution of $k$ TASEP, usually called basic coupling.

\subsection{HAD process}
Let $\Lambda:=\mathbb R/\mathbb Z$ be the one dimensional torus. The
Hammersley-Aldous-Diaconis (HAD) process (after \cite{H} and
\cite{AD}) is a stochastic evolution on finite subsets of $\Lambda$.
Let
$$
\Omega_N:=\left\{\underline x:=\{x_1,\dots ,x_N\}\ :\
x_i\in\Lambda;\ x_i\neq x_j \  \ i\neq j\right\}
$$
be the collection of all finite subsets of $\Lambda$ with $N$ points
and let $\Omega=\cup_N\Omega_N$. Labels to points are given in such
a way that $x_{i+1}$ is the nearest point of $\underline x$ to the
right of $x_i$. Given $\underline x\in \Omega$ an initial condition,
the HAD process preserves the number of points and is defined as
follows. Every point $x_i$ waits an exponential time of rate
$|(x_{i},x_{i+1}]|$ and then jumps to a point uniformly chosen in
$(x_i,x_{i+1}]$. This dynamics can be easily summarized from the
following generator
$$
Lf(\underline x)=\int_{\Lambda} du \ \left(f(\underline
x^u)-f(\underline x)\right)\ ,
$$
where the set $\underline x^u$ is defined from
\begin{equation}
x^u_i:=\left\{
\begin{array}{cc}
x_i & if \ u\notin (x_i,x_{i+1}]\ , \\
u & if \ u \in (x_i,x_{i+1}]\ .\\
\end{array}
\right. \nonumber
\end{equation}
This formula holds for any $i=1, \dots ,|\underline x |$, with the
convention $x_{|\underline x |+1}:=x_1$. With probability one all
the points stay distinct along the evolution.

The multiclass HAD process has not a simple and intuitive behavior.
As in the case of the TASEP a natural way to define it is through
the basic coupling. The state space is $I^{\uparrow}_{N_1,\dots
,N_k}$ defined as
$$
I^{\uparrow}_{N_1,\dots ,N_k}:=\left\{\left(\underline x^{(1)},\dots
,\underline x^{(k)}\right): \ \underline x^{(i)}\in \Omega_{N_i}; \
\underline x^{(i)}\subseteq \underline x^{(i+1)}\right\}\ .
$$
Here $0\leq N_1\leq N_2 \cdots \leq N_k$ are $k$ natural numbers.
This is a natural set to describe points in $\Lambda$ with
associated an integer class from $1$ to $k$. Points are the elements
of $\underline x^{(k)}$. The class associated to $ x^{(k)}_i$ is
$\inf \left\{j:\ x^{(k)}_i\in \underline x^{(j)}\right\}$. The
multiclass dynamics is defined from the following generator
$$
Lf(\underline x^{(1)},\dots ,\underline x^{(k)})=\int_{\Lambda} du \
\left(f(\underline x^{(i),u},\dots ,\underline
x^{(k),u})-f(\underline x^{(1)},\dots ,\underline x^{(k)})\right)
$$
that describes also a joint, inclusion preserving, evolution of $k$
HAD processes. This joint evolution is usually called basic
coupling.

Both for the TASEP and the HAD process a natural way to introduce
basic coupling is via a graphical construction where all the coupled
processes evolve using the same random marks. In the case of the HAD
process for example these marks are points of a rate one Poisson
point process on the cylinder $\Lambda\times \mathbb R^+$. See
\cite{FM1}, \cite{FM2}, \cite{FM3} for a more detailed description.

\section{Collapsing particles}
\label{colla}

To describe the invariant measures of the multiclass processes
previously introduced we have to explain a \emph{collapsing
procedure} introduced in \cite{A}. We start describing its action on
configurations of the TASEP.

Let us call
$$
X_{N,M}:=\left\{\eta\in X_N : \sum_{x\in \mathbb
Z_N}\eta(x)=M\right\}\ .
$$
Given $0\leq M_1\leq M_2\leq N$ two natural numbers, we define a
\emph{collapsing operator}
$$
\mathbb C: X_{N,M_1}\times X_{N,M_2}\to I_N^{2,\uparrow}
$$
that maps the pair of configurations $(\eta_1,\eta_2)$ into the pair
of configurations
$$
\mathbb C (\eta_1,\eta_2):=\left(C_{\eta_2}[\eta_1],\eta_2\right)\ .
$$
The collapsed configuration $C_{\eta_2}[\eta_1]$ is obtained from a
mass preserving (i.e. number of particle preserving) transformation
of the configuration $\eta_1$. The transformation is defined
algorithmically as follows. Give any order to the particles of the
configuration $\eta_1$. Move the first particle of this
configuration that is on a site, say $x$, such that $\eta_2(x)=0$,
to the first site on the right, say $y$, such that: $\eta_2(y)=1$
and $\eta_1(y)=0$. Update the $\eta_1$ configuration according to
this movement and iterate the procedure using at every step the same
order fixed at the beginning. The final configuration does not
depend on the specific order chosen.

We briefly discuss some of the properties of this algorithmic
transformation. For more details we refer the reader to \cite{A}
where this construction was introduced, and to \cite{FM1},
\cite{FM2}, \cite{FM3} for an interpretation in terms of queue
theory.

By definition we have that $C_{\eta_2}[\eta_1]\preceq \eta_2$. Let
$E(y,x)$ be the excess of $\eta_1$ particles in $[y,x]$ defined as
follows
\begin{equation}
E(y,x):=\sum_{z\in [y, x]}\left(\eta_1(z)-\eta_2(z)\right)\ .
%\label{flux}
\end{equation}
\begin{Le}
There is a positive flux of particles across the bond $e_x$ if and
only if there exists $y\in \mathbb Z_N$ such that $E(y,x)>0$.
\end{Le}
\begin{proof}
If $y$ is such that $E(y,x)>0$ then in $[y,x]$ there are more
particles of $\eta_1$ type than of $\eta_2$. During the collapsing
procedure some particles will necessarily flow out of $[y,x]$ and
this can happens only through $e_x$. Conversely let us suppose that
for any $y\in \mathbb Z_N$ we have $E(y,x)\leq 0$. Let us order the
particles of $\eta_1$ from right to left starting from x and let
$z_i$ be the site corresponding to the $i$ particle. Due to the fact
that $ E(z_1, x)\leq 0$, particle number $1$ will be allocated on a
site belonging to $[z_1,x]$. Due to the fact that $E(z_2, x)\leq 0$,
particle number $2$ will be allocated on a site belonging to
$[z_2,x]$ (different from the one of particle number $1$) and so on.
No particles will flow across $e_x$.
\end{proof}

Let $\left[\ \cdot \ \right]_+$ be the positive part defined as
\begin{equation}
[x]_+:=\left\{
\begin{array}{lc}
x & if \ x\geq 0\ ,\\
0 & otherwise \ .\\
\end{array}
\right.\nonumber
\end{equation}
\begin{Le}
The  total flux of particles across $e_x$ is
$$
J(x):=\sup_{y\in \mathbb Z_N}\left[E(y,x)\right]_+\ .
$$
\end{Le}
\begin{proof}
When $J(x)=0$ this follows directly from the previous lemma. If
$J(x)>0$ we can argue as follows. Clearly $J(x)$ is a lower bound of
the total flux, because for any $y$ the excess of $\eta_1$ particles
$E(y,x)$, if positive, has necessarily to flow across $e_x$. Let
$y^*$ be the first element of $\mathbb Z_N$ to the left of $x$ such
that $E(y^*,x)=J(x)>0$. There is no flux of particles across
$e_{y^*-1}$. This follows from the previous lemma and the fact that
for any $z$ it holds $E(z,y^*-1)\leq 0$. In fact, if $z\in
[x+1,y^*-1]$ then
$$
E(z,x)\leq J(x)=E(y^*,x)\ ,
$$
that implies $E(z,y^*-1)=E(z,x)-E(y^*,x) \leq 0$. If instead $z\in
[y^*,x]$ we have
\begin{eqnarray*}
& E(z,y^*-1)+E(y^*,x) = E(z,x)+E(x+1,x) \\
& \leq J(x)+M_1-M_2  =  E(y^*,x)+M_1-M_2\ ,
\end{eqnarray*}
that implies $E(z,y^*-1)\leq M_1-M_2\leq 0$. We established that all
the particles flowing across $e_x$ were originally in $[y^*,x]$.
From the characterization of $y^*$ we deduce immediately that
$E(y^*,z)>0$ for any $z\in [y^*,x]$. Let us order the particles of
$\eta_2$ type contained in $[y^*,x]$ from left to right and let
$z_i$ be the site corresponding to the particle number $i$. Remember
that there is an excess $J(x)$ of $\eta_1$ particles in this
interval. Due to the fact that $E(y^*,z_1)>0$ a particle of type
$\eta_1$ will be allocated in $z_1$. Due to the fact that
$E(y^*,z_2)>0$ a particle of $\eta_1$ type will be allocated in
$z_2$ and so on. At the end all the sites $z_i$ will be occupied by
$\eta_1$ particles and exactly the excess of particles $J(x)$ will
flow trough $e_x$.
\end{proof}

\begin{Le}
For any interval $[a,b]$ it holds
\begin{equation}
\sum_{x\in [a,b]}C_{\eta_2}[\eta_1](x)=\sum_{x\in [a,b]}
\eta_1(x)+J(a-1)-J(b)\ . \label{central}
\end{equation}
\end{Le}
\begin{proof}
This property follows directly from the conservation of mass.
Equation \eqref{central} simply states that the number of $\eta_1$
type particles that are at the end of the collapsing procedure in
the interval $[a,b]$ is obtained from the number of particles
present initially plus the number of particles entered from the left
side minus the number of particles exit from the right side.
\end{proof}

The collapsing procedure is defined in a similar way for
configurations of the HAD process. Let $0\leq N_1\leq N_2$ be two
integer numbers. We define the collapsing operator
$$
\mathbb C: \Omega_{N_1}\times \Omega_{N_2}\to I_{N_1,N_2}^{\uparrow}
$$
that maps the pair of configurations $(\underline x,\underline y)$
into the pair of configurations
$$
\mathbb C (\underline{x},\underline{y}):=\left(C_{\underline
y}[\underline x],\underline y\right)\ .
$$
The collapsed configuration $C_{\underline y}[\underline x]$ is
obtained moving to the right points of $\underline x$. The
transformation is defined algorithmically as follows. Give any order
to points of $\underline x$. Move the first point of $\underline x$
that does not belong to $\underline y$ to the nearest point of
$\underline y$ to the right that does not belong to $\underline x$.
Update the $\underline x$ configuration according to the previous
transformation and iterate the procedure. The final configuration
does not depend on the specific order chosen.

Also in this case it is possible to define an excess $E(u,v)$ of $
\underline{x}$ points in the interval $[u,v]$
$$
E(u,v):=\Big|\left\{\underline{x}\cap[u,v]\right\}\Big|-\Big|\left\{\underline{y}\cap
[u,v]\right\}\Big|
$$
and consequently a flux of $ \underline{x}$ particles at $v$
$$
J(v):=\sup_{u\in \Lambda}\left[E(u,v)\right]_+ \ .
$$
Note that $J$ is right continuous i.e. it holds
$J(v)=J(v^+):=lim_{\epsilon\downarrow 0}\ J(v+\epsilon)$. All the
lemmas previously listed, appropriately reformulated, holds also in
this case. We do not go into details here because the collapsing
procedure will be generalized to a wider framework in section
\ref{collames}. We only write down the analogous of equation
\eqref{central} in this case
\begin{equation}
\Big|\left\{C_{\underline
y}[\underline{x}]\cap[u,v]\right\}\Big|=\Big|\left\{\underline{x}\cap[u,v]\right\}\Big|+J(u^-)-J(v)\
, \label{central2}
\end{equation}
where $J(u^-)$ is the left limit of $J$ at $u$, i.e.
$J(u^-):=lim_{\epsilon\downarrow 0}\ J(u-\epsilon)$.

\section{Invariant measures}
\label{invcolla}

\subsection{TASEP} For the TASEP the number of particles
$\sum_{x\in \mathbb Z_N}\eta(x)$ is a conserved quantity. For any
fixed integer $M\leq N$, the TASEP with $M$ particles is an
irreducible finite state Markov chain on $X_{N,M}$ and consequently
has a unique invariant measure. The process has then a one parameter
(i.e. $M$) family of invariant measures that is easily seen to
coincide with the family of uniform distributions of $M$ particles
on $\mathbb Z_N$
\begin{equation}
\nu_N^M(\eta)=\left\{
\begin{array}{cc}
\binom{N}{M}^{-1} & if \ \ \eta\in X_{N,M}\ ,\\
0& otherwise\ .\\
\end{array}
\right. \nonumber
\end{equation}
These are all the extremal invariant measures; all the remaining are
obtained as convex combinations.

The 2-class TASEP conserves the number of first class particles
$\sum_{x\in \mathbb Z_N}\eta_1(x)$ and the number of second class
particles $ \sum_{x\in \mathbb Z_N}(\eta_2(x)-\eta_1(x))$. For any
fixed pair of non negative integer numbers $\Delta_1$ and $\Delta_2$
such that $M_1:=\Delta_1$ and $M_2:=\Delta_1+\Delta_2\leq N$, the
TASEP with $\Delta_1$ first class particles and $\Delta_2$ second
class particles is an irreducible finite state Markov chain and has
an unique invariant measure.

The result in \cite{A} states that this invariant measure is
\begin{equation}
\left(\nu_N^{M_1}\times \nu_N^{M_2}\right)\circ \mathbb C^{-1}\ .
\label{cacca}
\end{equation}

We use the symbol $\times$ to indicate the product of measures. In
general given a measure $\mu$ and a measurable map $T$ with the
symbol $\mu\circ T^{-1}$ we denote the pull-back measure defined
from
$$
\left(\mu\circ T^{-1}\right)(A):=\mu\left(T^{-1}(A)\right)
$$
for any measurable set $A$.

We remark that \eqref{cacca} is a measure on $I^{2,\uparrow}_N$.
This two parameters ($M_1$ and $M_2$) family of invariant measures
constitutes all the extremal invariant measures.

To give a combinatorial representation of the invariant measures of
the $k$-class TASEP we need to extend the collapsing procedure of
section \ref{colla}. This extension to the case of more than two
classes of particles is contained in \cite{FM1} and further
discussed in \cite{FM2}, \cite{FM3}. Let $\Delta_1, \dots ,\Delta_k$
be $k$ non negative integer numbers such that
$\sum_{i=1}^k\Delta_i\leq N$. Let us call also
$M_j:=\sum_{i=1}^j\Delta_i$. We define a collapsing operator
$$
\mathbb C_k: X_{N,M_1}\times X_{N,M_2}\times \cdots \times
X_{N,M_k}\to I_N^{k,\uparrow}
$$
that associates to the configurations $(\eta_1,\dots ,\eta_k)$ the
configurations
$$
(\xi_1,\dots ,\xi_k):=\mathbb C_k(\eta_1, \dots ,\eta_k)
$$
defined as follows. The configuration $\xi_k$ coincides with
$\eta_k$. The configuration $\xi_{k-1}$ coincides with
$C_{\eta_k}[\eta_{k-1}]$. The configuration $\xi_{k-2}$ coincides
with $C_{\eta_k}[C_{\eta_{k-1}}[\eta_{k-2}]]$. In general the
configuration $\xi_{k-j}$ is obtained from the composition of $j$
collapsing procedures
\begin{equation}
\xi_{k-j}:=C_{\eta_k}\left[C_{\eta_{k-1}}\left[\dots
C_{\eta_{k-j+1}}\left[\eta_{k-j}\right]\dots\right]\right]\ .
\label{cat1}
\end{equation}
Obviously according to this definition $\mathbb C=\mathbb C_2$. The
result contained in \cite{FM1} states that the invariant measure for
the TASEP with $\Delta_i$ i-class particles $(i=1,\dots ,k$) is
$$
\left(\nu_N^{M_1}\times \nu_N^{M_2}\times \cdots \times
\nu_N^{M_k}\right)\circ \mathbb C_k^{-1}\ .
$$
This k-parameter family of invariant measures constitutes all the
extremal invariant measures.

\subsection{HAD}
The HAD process conserves the number of points. In $\Omega_N$ the
unique invariant measure $\mu_N$ is given by the support of the
values of $N$ i.i.d random variables uniform in $\Lambda$.
Equivalently the invariant measure is a uniform Poisson point
process in $\Lambda$, conditioned to have $N$ points. These are all
the extremal invariant measures, all the remaining are obtained as
convex combinations.

The $2$-class HAD process conserves the number of first class points
and the number of second class points. Let $\Delta_1$ and $\Delta_2$
be two non negative integer numbers and let $N_1:=\Delta_1$ and
$N_2:=\Delta_1+\Delta_2$. The $2$-class HAD process with $\Delta_1$
points of first class and $\Delta_2$ points of second class has a
unique invariant measure that has a combinatorial representation in
terms of the collapsing operator. In fact the result in \cite{FM3}
states that this invariant measure is
$$
\left(\mu_{N_1}\times \mu_{N_2}\right)\circ \mathbb C^{-1}\ .
$$
These are all the extremal invariant measures, all the remaining are
obtained as convex combinations.

The generalization of this combinatorial construction to the case of
more than two classes of points is described also in \cite{FM3}. We
proceed as in the case of the TASEP. Let $\Delta_1, \Delta_2,\dots
\Delta_k$ be non negative integer numbers and call
$N_j:=\sum_{i=1}^j\Delta_i$. We define a collapsing operator
$$
\mathbb C_k: \Omega_{N_1}\times \Omega_{N_2 }\times \cdots \times
\Omega_{N_k}\to I_{N_1,\dots ,N_k}^{\uparrow}
$$
that associates to the configurations $\left(\underline
x^{(1)},\dots ,\underline x^{(k)}\right)$ the configurations
$$
\left(\underline y^{(1)},\dots ,\underline y^{(k)}\right):=\mathbb
C_k\left(\underline x^{(1)},\dots ,\underline x^{(k)}\right)
$$
defined as follows. The configuration $\underline y^{(k)}$ coincides
with $\underline x^{(k)}$. The configuration $\underline y^{(k-1)}$
coincides with $C_{\underline x^{(k)}}[\underline x^{(k-1)}]$. The
configuration $\underline y^{(k-2)}$ coincides with $C_{\underline
x^{(k)}}[C_{\underline x^{(k-1)}}[\underline x^{(k-2)}]]$.

In general the configuration $\underline{y}^{(k-j)}$ is obtained
from the composition of $j$ collapsing procedures
\begin{equation}
\underline{y}^{(k-j)}:=C_{\underline{x}^{(k)}}\left[C_{\underline{x}^{(k-1)}}\left[\dots
C_{\underline{x}^{(k-j+1)}}\left[\underline{x}^{(k-j)}\right]\dots\right]\right]\
. \label{cat2}
\end{equation}
The result in \cite{FM3} states that the invariant measure for the
HAD process with $\Delta_i$ i-class particles $(i=1,\dots ,k$) is
$$
\left(\mu_{N_1}\times \mu_{N_2}\times \cdots \times
\mu_{N_k}\right)\circ \mathbb C_k^{-1}\ .
$$
These are all the extremal invariant measures, all the remaining are
obtained as convex combinations.

\section{Collapsing measures}
\label{collames} We start with some definitions. The set of positive
measures on $\Lambda$ will be denoted as $\mathcal
 M$. With $\mathcal M^0$ we will denote the subset of $\mathcal M$
of measures $\rho$ absolutely continuous with respect to Lebesgue
measure and with $\mathcal M^{0,b}$ the subset of $\mathcal M^0$
containing the elements such that their densities satisfy the
condition
$$
0\leq \frac{d\rho}{du}\leq 1 \ \ \ \ \ \ \ a.e.\ .
$$
With abuse of notation we will indicate with $\rho$ both a generic
element of $\mathcal M^0$ and its density. Finally, given $m\in
\mathbb R^+$, we call $\mathcal M_m:=\left\{\rho\in \mathcal M : \
\int_{\Lambda}\ d\rho=m\right\}$. Likewise we set $\mathcal
M^0_m:=\mathcal M^0\bigcap \mathcal M_m$ and $\mathcal
M^{0,b}_m:=\mathcal M^{0,b}\bigcap \mathcal M_m$.

We define a partial order $\preceq$ on $\mathcal M$ saying that
$$
\rho_1\preceq \rho_2 \ \  \Leftrightarrow \ \ \int_A d\rho_1\leq
\int_A d\rho_2
$$
for any measurable $A\subseteq \Lambda$. We then
call
$$
I^{k,\uparrow}:=\left\{(\rho_1,\dots ,\rho_k): \rho_i\in \mathcal M;
\ \rho_i\preceq \rho_{i+1}\right\}\ .
$$

We want to generalize to this framework the algorithmic
constructions illustrated in section \ref{colla}. We are dealing no
more with configurations of particles, but with positive measures.
Given $\rho_1\in \mathcal M_{m_1}$ and $\rho_2\in \mathcal M_{m_2}$
with $m_1\leq m_2$ we want to define a collapsing operator
$$
\mathbb C:\mathcal M_{m_1}\times \mathcal M_{m_2}\to I^{2,\uparrow}
$$
that associates to the pair $(\rho_1,\rho_2)$ the pair
$$
\mathbb C(\rho_1,\rho_2):=\left(C_{\rho_2}[\rho_1],\rho_2\right)\ ,
$$
where the collapsed measure $C_{\rho_2}[\rho_1]\preceq \rho_2$ is
obtained \emph{moving mass of $\rho_1$ to the right}. The natural
way to define such a procedure is through a generalization of
equations \eqref{central} and \eqref{central2}. We start defining
the excess of mass $E(u,v)$ of the measure $\rho_1$ in the interval
$[u,v]$ as
$$
E(u,v):=\int_{[u,v]}d\rho_1-\int_{[u,v]}d\rho_2\ .
$$
By definition $E(u,v)$ is right continuous in $v$ (i.e.
$E(u,v^+):=\lim_{\epsilon \downarrow 0}E(u,v+\epsilon)=E(u,v)$) and
left continuous in $u$ (i.e. $E(u^-,v)=E(u,v)$) . It satisfies also
some simple addition rules
\begin{eqnarray*}
&E(a,c)=E(a,b^-)+E(b,c) \ \ \ \ \ \ \forall \ b\in [a,c]\ ,\\
&E(a,c)=E(a,b)+E(b^+,c) \ \ \ \ \ \  \forall \ b\in [a,c]\ .
\end{eqnarray*}
In the above formulas we use the convention $E(u,u^-)=E(u^+,u):=0$.
Then we introduce the flux of mass across $v\in \Lambda$ defined as
$$
J(v):=\sup_{u\in \Lambda}\left[E(u,v)\right]_+ \ .
$$
We give the following definition of the collapsing operator.
\begin{Def}
\label{deftot} The collapsed measure $C_{\rho_2}[\rho_1]$ is defined
in such a way that for any interval $(a,b]$ it holds
\begin{equation}
\int_{(a,b]}dC_{\rho_2}[\rho_1]=\int_{(a,b]}d\rho_1+J(a)-J(b)\ .
\label{def}
\end{equation}
\end{Def}
Note that the action on intervals of this type completely defines
the measure.

Let
\begin{equation}
\mathcal J:=\left\{v\in \Lambda : \ \exists u\in \Lambda \ s.t. \
E(u,v)>0\right\}\ . \label{dafare}
\end{equation}
Given $v\in \mathcal J$, let $u$ such that $E(u,v)>0$. From the
right continuity in $v$ we have $E(u,v^+)=E(u,v)>0$ and consequently
there exists an $\epsilon >0$ such that $v+\delta \in \mathcal J$
for any $0\leq \delta \leq \epsilon$. From this we can deduce that
$\mathcal J=\cup_i\mathcal J_i$. Where $\mathcal J_i$ are at most
countable many disjoint intervals either of the type $[l_i,r_i)$ or
$(l_i,r_i)$. This fact can be proved with an argument very similar
to the one used to characterize open sets on $\mathbb R$ (see for
example section II.11 of \cite{KF}). The condition $m_1\leq m_2$
implies the fact that the strict inclusion $\mathcal J\subset
\Lambda$ holds. A sketch of the proof is as follows. Let, by
contradiction, assume that $\mathcal J=\Lambda$. Then for any $v\in
\Lambda$ you can prove there exist an interval $\mathcal U_v$ whose
right closed boundary is $v$, it is either open or closed at the
left boundary and is such that
$$
0<J(v)=\int_{\mathcal U_v}d\rho_1-\int_{\mathcal U_v}d\rho_2\ .
$$
Moreover any other interval with the same property is contained
inside $\mathcal U_v$. Then for any $v_1\neq v_2$ we have that
either $\mathcal U_{v_1}\cap \mathcal U_{v_2}=\emptyset$ or they are
one contained inside the other. We define $\mathcal
W_v:\cup_{\left\{w\in \Lambda: \ v\in \mathcal U_w\right\}}\mathcal
U_w$. Given $v_1\neq v_2$ then either $\mathcal W_{v_1}=\mathcal
W_{v_2}$ or they are disjoint. Moreover
$$
\int_{\mathcal W_v}d\rho_1-\int_{\mathcal W_v}d\rho_2>0 \ \ \ \
\forall v\in \Lambda
$$
so that at most countable different $\mathcal W_v$ can exist and
they form a partition of $\Lambda$. We finally have
\begin{eqnarray*}
& 0\geq m_1-m_2=\int_{\Lambda}d\rho_1-\int_{\Lambda}d\rho_2 \\
&=\sum_{v_i}\left(\int_{\mathcal W_{v_i}}d\rho_1-\int_{\mathcal
W_{v_i}}d\rho_2\right)>0\ ,
\end{eqnarray*}
a contradiction. An alternative route is obtained showing with
arguments similar to the ones of the next lemma that there exists an
element of $\Lambda$ in a neighborhood of the left boundary of any
$\mathcal U_v$ that does not belong to $\mathcal J$.

\begin{Le}
\label{stralemma} Consider $v\in \mathcal J_i$. If $\mathcal
J_i=[l_i,r_i)$ then
$$
J(v)=E(l_i,v)\ .
$$
If instead $\mathcal J_i=(l_i,r_i)$ then
$$
J(v)=E(l_i^+,v)\ .
$$
\end{Le}
\begin{proof}
Here and hereafter in some proofs we need to distinguish the two
cases when $\mathcal J_i=[l_i,r_i)$ or $\mathcal J_i=(l_i,r_i)$. We
will give the proofs in the case $\mathcal J_i=[l_i,r_i)$. The
proofs for the other case are analogous. Given $v\in \mathcal J_i$
we consider $w_n$ a maximizing sequence for $J(v)$, i.e. a sequence
such that
$$
\lim_{n\to \infty}E(w_n,v)=J(v)\ .
$$
If $w_n\in (v,l_i)$ then
\begin{equation}
E(w_n,v)=E(w_n,l_i^-)+E(l_i,v)\leq E(l_i,v)\ , \label{summa}
\end{equation}
where we used the definition of $\mathcal J$ and the fact that there
exists an element of $\mathcal J^c$, the complement of $\mathcal J$
in $\Lambda$, in any neighborhood of $l_i$. Without loss of
generality we can then consider $w_n\in [l_i,v]$. Note also that the
simple argument in \eqref{summa} also implies that for any $v\in
\mathcal J_i$ there exists $u\in[l_i,v]$ such that $E(u,v)>0$.

We prove now that for any $w\in \mathcal J_i$ it holds $E(l_i,w)\geq
0$. Let us introduce
$$
\mathcal N:=\left\{w\in \mathcal J_i: \ E(l_i,w)<0\right\}\ .
$$
Due to the fact that $E(u,v)$ is right continuous in $v$ we have
that $\mathcal N$ is the union of at most countable many intervals
either of the type $[a_j,b_j)$ or of the type $(a_j,b_j)$. We give
the proof in the case there are only a finite number of intervals
and $a_1$, the nearest boundary element to the right of $l_i$, is
the boundary  of an interval of the type $(a_1,b_1)$. As before, the
proof for the remaining cases is analogous. As already showed there
exists $u\in [l_i,a_1]$ with $E(u,a_1)>0$. Then we have
$$
0\geq E(l_i,a_1^+)=E(l_i,a_1)=E(l_i,u^-)+E(u,a_1)\geq E(u,a_1)\ ,
$$
a contradiction. This imply that $\mathcal N=\emptyset$. Finally we
get for the maximizing sequence $w_n\in [l_i,v]$ that
$$
E(l_i,v)=E(l_i,w_n^-)+E(w_n,v)\geq E(w_n,v)\ .
$$
This means that $\widetilde{w}_n:=l_i$ is a maximizing sequence and
this implies the first statement of the lemma.
\end{proof}
Note that a direct consequence of this lemma is that $J$ is right
continuous, an important fact to have that definition \ref{deftot}
is well posed.

It is possible to introduce a measure $\gamma$ defined from
$$
\int_{(a,b]}d\gamma:=J(b)-J(a)\ .
$$
The measure $\gamma$ is not a positive measure. In fact we have
\begin{equation}
\int_{\Lambda}d\gamma=\lim_{\epsilon \downarrow
o}\int_{(v+\epsilon,v]}d\gamma=\lim_{\epsilon \downarrow
0}\left(J(v)-J(v+\epsilon)\right)=0\ . \label{zero}
\end{equation}

Definition \eqref{def} then becomes
$$
C_{\rho_2}[\rho_1]=\rho_1-\gamma
$$
and using \eqref{zero} we derive the conservation of mass
$$
\int_{\Lambda}
dC_{\rho_2}[\rho_1]=\int_{\Lambda}d\rho_1-\int_{\Lambda}d\gamma=\int_{\Lambda}d\rho_1\
.
$$
We give now a simple representation of the collapsed measure
$C_{\rho_2}[\rho_1]$. Similar representations hold also for
collapsed configurations of particles in the TASEP and points in the
HAD process but their generalization to the case of positive
measures is not straightforward.

\begin{Le}
\label{lemme}
The collapsed measure $C_{\rho_2}[\rho_1]$ has the
following representation
\begin{equation}
C_{\rho_2}[\rho_1]=\rho_1\chi_{_{\mathcal
J^c}}+\rho_2\chi_{_{\mathcal J}}+\sum_i\left(\int_{\mathcal
J_i}d\rho_1-\int_{\mathcal J_i}d\rho_2\right)\delta_{r_i}\ ,
\label{semplice}
\end{equation}
where with the symbol $\chi_{_A}$ we denote the characteristic
function of the set $A$ and with $\delta_v$ the delta measure in
$v$.
\end{Le}
\begin{proof}
We show that the weight given to any interval $(a,b]$ from the
measure defined by the right hand side of \eqref{def} coincide with
the weight given to the same interval from the measure on the right
hand side of \eqref{semplice}. This implies that the two measures
coincide.

We need to verify the following identity
\begin{eqnarray*}
& &\int_{(a,b]}d\rho_1+J(a)-J(b)\\
&= &\int_{(a,b]\cap\mathcal J^c}d\rho_1+\int_{(a,b]\cap\mathcal
J}d\rho_2+\sum_{r_i\in(a,b]}\left(\int_{\mathcal
J_i}d\rho_1-\int_{\mathcal J_i}d\rho_2\right)\ .
\end{eqnarray*}
This is equivalent to
\begin{eqnarray}
& &\int_{(a,b]\cap \mathcal J}d\rho_1-\int_{(a,b]\cap \mathcal
J}d\rho_2=-\left(\int_{[l_{i^*},a]}d\rho_1-
\int_{[l_{i^*},a]}d\rho_2\right)\chi_{_\mathcal J}(a)\nonumber \\
&+ &\left(\int_{[l_{j^*},b]}d\rho_1-
\int_{[l_{j^*},b]}d\rho_2\right)\chi_{_\mathcal J}(b)
+\sum_{r_i\in(a,b]}\left(\int_{\mathcal J_i}d\rho_1-\int_{\mathcal
J_i}d\rho_2\right)\ . \label{nuovadim}
\end{eqnarray}
In the above formula, when $a\in \mathcal J$ we call $\mathcal
J_{i^*}$ the interval to which it belongs and we assume it is of the
type $[l_{i^*},r_{i^*})$, when $b\in \mathcal J$ we call $\mathcal
J_{j^*}$ the interval to which it belongs and we assume it is of the
type $[l_{j^*},r_{j^*})$. The proof in the remaining cases is
similar. Note that it is possible to have $i^*=j^*$.

The sum on the right hand side of \eqref{nuovadim} can be written in
the following way.
\begin{eqnarray}
& &\sum_{r_i\in(a,b]}\left(\int_{\mathcal J_i}d\rho_1-\int_{\mathcal
J_i}d\rho_2\right) \nonumber \\
&=&\sum_{\mathcal J_i\subseteq (a,b]}\left(\int_{\mathcal
J_i}d\rho_1-\int_{\mathcal J_i}d\rho_2\right)+\left(\int_{\mathcal
J_{i^*}}d\rho_1- \int_{\mathcal
J_{i^*}}d\rho_2\right)\chi_{_\mathcal J}(a)\ . \label{nuovadeco}
\end{eqnarray}
Using this formula the right hand side of \eqref{nuovadim} becomes
\begin{eqnarray*}
& &\sum_{\mathcal J_i\subseteq (a,b]}\left(\int_{\mathcal
J_i}d\rho_1-\int_{\mathcal
J_i}d\rho_2\right)+\left(\int_{(a,r_{i^*})}d\rho_1-
\int_{(a,r_{i^*})}d\rho_2\right)\chi_{_\mathcal J}(a)\\
&+&\left(\int_{[l_{j^*},b]}d\rho_1-
\int_{[l_{j^*},b]}d\rho_2\right)\chi_{_\mathcal J}(b)=
\sum_{i}\left(\int_{\mathcal J_i\cap (a,b]}d\rho_1-\int_{\mathcal
J_i\cap (a,b]}d\rho_2\right)\ ,
\end{eqnarray*}
that clearly coincides with the left hand side of \eqref{nuovadim}.
\end{proof}

We have the following simple properties.
\begin{Le}
The measure
$C_{\rho_2}[\rho_1]$ is a positive measure.
\end{Le}
\begin{proof}
We use the representation \eqref{semplice}. The first two terms on
the right hand side are clearly positive. In the case $\mathcal
J_i=[l_i,r_i)$ we have that
$$
\int_{\mathcal J_i}d\rho_1-\int_{\mathcal
J_i}d\rho_2=E(l_i,r_i^-)\geq 0 \ ,
$$
as shown during the proof of lemma \ref{stralemma}. This implies
that also the third term is positive. An analogous argument holds in
the case $\mathcal J_i=(l_i,r_i)$.
\end{proof}

\begin{Le}
It holds
\begin{equation}
C_{\rho_2}[\rho_1]\preceq \rho_2\ . \label{minormes}
\end{equation}
\end{Le}
\begin{proof}
Using formula \eqref{semplice} we obtain that condition
\eqref{minormes} is equivalent to
\begin{equation}
\rho_1\chi_{_{\mathcal J^c}}-\rho_2\chi_{_{\mathcal
J^c}}+\sum_i\left(\int_{\mathcal J_i}d\rho_1-\int_{\mathcal
J_i}d\rho_2\right)\delta_{r_i} \preceq 0 \ .\label{nuovadim2}
\end{equation}
We show that the weight associated to any interval $(a,b]$ from the
measure on the left hand side of \eqref{nuovadim2} is negative and
this implies the statement of the lemma.

We need to prove that for any $(a,b]$ it holds
\begin{equation}
\int_{(a,b]\cap \mathcal J^c}d\rho_1-\int_{(a,b]\cap \mathcal
J^c}d\rho_2+\sum_{r_i\in(a,b]}\left(\int_{\mathcal
J_i}d\rho_1-\int_{\mathcal J_i}d\rho_2\right)\leq 0 \ .
\label{nuovadim3}
\end{equation}
We use the fact that
\begin{eqnarray}
& &\int_{(a,b]\cap \mathcal J^c}d\rho_1-\int_{(a,b]\cap \mathcal
J^c}d\rho_2 \nonumber \\
&=&\int_{(a,b]}d\rho_1-\int_{(a,b]}d\rho_2-\sum_{i}\left(\int_{\mathcal
J_i\cap
(a,b]}d\rho_1-\int_{\mathcal J_i\cap (a,b]}d\rho_2\right)\nonumber \\
&=& \int_{(a,b]}d\rho_1-\int_{(a,b]}d\rho_2-\sum_{ \mathcal
J_i\subseteq (a,b]}\left(\int_{\mathcal J_i}d\rho_1-\int_{\mathcal J_i}
d\rho_2\right)\nonumber \\
&-&\left(\int_{(a,r_{i^*})}d\rho_1-
\int_{(a,r_{i^*})}d\rho_2\right)\chi_{_\mathcal
J}(a)-\left(\int_{[l_{j^*},b]}d\rho_1-
\int_{[l_{j^*},b]}d\rho_2\right)\chi_{_\mathcal J}(b)\ . \label{uhm}
\end{eqnarray}
In the above formulas we use for the labels $i^*$ and $j^*$ and the
corresponding intervals, the same convention as in lemma
\ref{lemme}. Using \eqref{uhm} and equation \eqref{nuovadeco} we
obtain that the left hand side of \eqref{nuovadim3} is equal to
\begin{eqnarray}
& &\int_{(a,b]}d\rho_1-\int_{(a,b]}d\rho_2+
\left(\int_{[l_{i^*},a]}d\rho_1-
\int_{[l_{i^*},a]}d\rho_2\right)\chi_{_\mathcal J}(a)\nonumber
\\
&-& \left(\int_{[l_{j^*},b]}d\rho_1-
\int_{[l_{j^*},b]}d\rho_2\right)\chi_{_\mathcal J}(b)\ .
\label{ramna}
\end{eqnarray}
In the case $i^*\neq j^*$ we have that \eqref{ramna} can be written
as
\begin{eqnarray*}
& & E(l_{i^*},l_{j^*}^-)\chi_{_\mathcal J}(a)\chi_{_\mathcal
J}(b)+E(l_{i^*},b)\chi_{_\mathcal J}(a)\chi_{_{\mathcal J^c}}(b)\\
& &+ E(a^+,l_{j^*}^-)\chi_{_{\mathcal J^c}}(a)\chi_{_\mathcal
J}(b)+E(a^+,b)\chi_{_{\mathcal J^c}}(a)\chi_{_{\mathcal J^c}}(b)
\end{eqnarray*}
and all the terms are non positive.

The case $i^*=j^*$ can happen only when $a\in \mathcal J$ and $b\in
\mathcal J$. In this case we have that \eqref{ramna} is equal either
to $0$ or to
$$
\int_{\Lambda}d\rho_1-\int_{\Lambda}d\rho_2\leq 0\ .
$$
\end{proof}

When both $\rho_1$ and $\rho_2$ belong to $\mathcal M^0$ then
$E(u,v)$ is continuous in $u$ and $v$ and consequently
$E(l_i,r_i^-)=E(l_i,r_i)$. From what we proved in lemma
\ref{stralemma} we have that $E(l_i,r_i^-)\geq 0$. From the fact
that $r_i\in \mathcal J^c$ we have that $E(l_i,r_i)\leq 0$. This
implies that $E(l_i,r_i)=0$. Then the third term on the right hand
side of \eqref{semplice} is not present and in this case we simply
have
$$
C_{\rho_2}[\rho_1]=\rho_1\chi_{_{\mathcal
J^c}}+\rho_2\chi_{_{\mathcal J}} \ .
$$
This means that also $C_{\rho_2}[\rho_1]\in \mathcal M^0$ with a
density a.e. given by
\begin{equation}
C_{\rho_2}[\rho_1](v)=\left\{
\begin{array}{ll}
\rho_2(v) & if \  \ v\in \mathcal J \ ,\\
\rho_1(v) & if \ \ v\in \mathcal J^c \ .\\
\end{array}
\right. \label{azasscont}
\end{equation}

\vskip0.5cm

The collapsing operator $\mathbb C$ is not continuous with respect
to the weak topology. This can be easily seen from the following
example. Consider the sequences of measures
$\rho_1^{(n)}:=\delta_{\left(\frac{1}{2}+\frac{1}{n}\right)}$ and
$\rho_2^{(n)}:=\delta_{\frac{1}{2}}+\delta_{\frac{3}{4}}$. Clearly
we have
$$
C_{\rho_2^{(n)}}[\rho_1^{(n)}]=\delta_{\frac{3}{4}}\ .
$$
Moreover it holds
$$
\rho_1^{(n)}\stackrel{n\to \infty}{\to} \rho_1 \ \ \ \ \ \  and \ \
\ \ \ \ \rho_2^{(n)}\stackrel{n\to \infty}{\to} \rho_2
$$
with $\rho_1:=\delta_{\frac{1}{2}}$,
$\rho_2:=\delta_{\frac{1}{2}}+\delta_{\frac{3}{4}}$ and convergence
is in the weak topology. We have also
$$
C_{\rho_2}[\rho_1]=\delta_{\frac{1}{2}}\ .
$$
As a consequence
$$
\lim_{n\to \infty}\mathbb C(\rho_1^{(n)},\rho_2^{(n)})\neq \mathbb
C(\rho_1,\rho_2)=\mathbb C\left(\lim_{n \to
\infty}(\rho_2^{(n)},\rho_1^{(n)})\right)\ .
$$
Nevertheless it holds the following continuity result.
\begin{Le}
\label{lemmino} The collapsing operator $\mathbb C$ is continuous
with respect to the weak topology on $\mathcal M^0_{m_1}\times
\mathcal M^0_{m_2}$.
\end{Le}
\begin{proof}
Let $(\rho_1^{(n)},\rho_2^{(n)})\in \mathcal M_{m_1}\times \mathcal
M_{m_2}$ be a sequence of measures weakly converging to
$(\rho_1,\rho_2)\in \mathcal M^0_{m_1}\times \mathcal M^0_{m_2}$.
Then for any fixed $v\in \Lambda$ we define the following
nondecreasing continuous functions on the interval $[0,1]$
\begin{equation}
G_1(u):=\left\{
\begin{array}{cc}
\int_{[v-u,v]}d\rho_1 & if \ u\in [0,1)\ ,\\
m_1 & if \ u=1\ ,
\end{array}
\right.\nonumber
\end{equation}
\begin{equation}
G_2(u):=\left\{
\begin{array}{cc}
\int_{[v-u,v]}d\rho_2 & if \ u\in [0,1)\ ,\\
m_2 & if \ u=1\ .
\end{array}
\right.\nonumber
\end{equation}
Likewise using the measures $\rho_1^{(n)}$ and $\rho_2^{(n)}$, we
define the functions $G^{(n)}_1$ and $G^{(n)}_2$ that are
nondecreasing but not necessarily continuous. The weak convergence
implies the pointwise convergence of $G^{(n)}_1$ to $G_1$ and
$G^{(n)}_2$ to $G_2$. Monotonicity of all the $G$, and continuity of
the limit functions imply also the uniform convergence. Consider for
example $G^{(n)}_1$ and $G_1$. Fix an arbitrary $\epsilon$, let
$$
u_i:=\inf \left\{u: \ G_1(u)=i\epsilon \right\}; \ \ \ \ \ i=0,\dots
,\left[\frac{m_1}{\epsilon}\right]
$$
and define also $u_{\left[\frac{m_1}{\epsilon}\right]+1}:=1$. For
any $i$ let $n(i)$ be such that for any $n>n(i)$ it holds
$|G_1^{(n)}(u_i)-G_1(u_i)|\leq \epsilon$. Let also
$n^*:=\max_i\left\{n(i)\right\}.$ Monotonicity says that when $u\in
[u_i,u_{i+1})$
$$
G_1^{(n)}(u_i)\leq G_1^{(n)}(u) \leq G_1^{(n)}(u_{i+1})
$$
and
$$
G_1(u_i)\leq G_1(u) \leq G_1(u_{i+1})\ .
$$
These inequalities imply
$$
|G_1^{(n)}(u)-G_1(u)|\leq \max \left\{
|G_1^{(n)}(u_i)-G_1(u_{i+1})|,
|G_1^{(n)}(u_{i+1})-G_1(u_i)|\right\}\ .
$$
When $n>n^*$ both terms inside the $\max$ are $\leq 2\epsilon$. Let
us show this for example for the second one. It holds
$$
|G_1^{(n)}(u_{i+1})-G_1(u_i)|\leq
|G_1^{(n)}(u_{i+1})-G_1(u_{i+1})|+|G_1(u_{i+1})-G_1(u_i)|
$$
and both terms on the right hand side are $\leq \epsilon$. This
implies the uniform convergence of $G^{(n)}_1$ to $G_1$. The uniform
convergence of $G^{(n)}_1$ and $G^{(n)}_2$ implies the uniform
convergence of $\left[G_1^{(n)}-G_2^{(n)}\right]_+$ to
$\left[G_1-G_2\right]_+$. This implies the convergence of
\begin{equation}
J^{(n)}(v)=\sup_{u\in[0,1]}\left[G_1^{(n)}(u)-G_2^{(n)}(u)\right]_+
\label{el}
\end{equation}
to
\begin{equation}
J(v)=\sup_{u\in[0,1]}\left[G_1(u)-G_2(u)\right]_+ \label{sad}
\end{equation}
for any $v\in \Lambda$. To prove this simple fact we see that,
calling $u^*$ a maximum point in \eqref{sad},
\begin{eqnarray*}
& \liminf_{n\to \infty}J^{(n)}(v)=\liminf_{n\to
\infty}\left(\sup_{u\in[0,1]}\left[G_1^{(n)}(u)-G_2^{(n)}(u)\right]_+\right)\\
&\geq\liminf_{n\to
\infty}\left[G_1^{(n)}(u^*)-G_2^{(n)}(u^*)\right]_+=J(v)\ .
\end{eqnarray*}
On the opposite direction, let us call $w_k^n$ a maximizing sequence
in \eqref{el}, i.e. a sequence such that
$$
J^{(n)}(v)=\lim_{k\to
\infty}\left[G_1^{(n)}(w^n_k)-G_2^{(n)}(w^n_k)\right]_+\ .
$$
Consider now a sequence $k(n)$ such that
$$
\left|J^{(n)}(v)-\left[G_1^{(n)}(w^n_{k(n)})-G_2^{(n)}(w^n_{k(n)})\right]_+\right|<\frac{1}{n}\
.
$$
It holds
$$
\limsup_{n\to \infty}J^{(n)}(v)=\limsup_{n\to
\infty}\left[G_1^{(n)}(w^n_{k(n)})-G_2^{(n)}(w^n_{k(n)})\right]_+\ .
$$
We then obtain
\begin{eqnarray*}
& \limsup_{n\to \infty}J^{(n)}(v)=\limsup_{n\to
\infty}\left[G_1^{(n)}(w^n_{k(n)})-G_2^{(n)}(w^n_{k(n)})\right]_+\\
&=
\limsup_{n\to\infty}\left[G_1(w^n_{k(n)})-G_2(w^n_{k(n)})\right]_+\leq
J(v)\ .
\end{eqnarray*}
We conclude the proof showing that
\begin{eqnarray*}
&\lim_{n\to
\infty}\int_{(a,b]}dC_{\rho_2^{(n)}}\left[\rho^{(n)}_1\right]\\
&=\lim_{n\to
\infty}\left(\int_{(a,b]}d\rho_1^{(n)}+J^{(n)}(a)-J^{(n)}(b)\right)\\
&=\int_{(a,b]}d\rho_1+J(a)-J(b)=\int_{(a,b]}dC_{\rho_2}\left[\rho_1\right]\
.
\end{eqnarray*}
The convergence of $\int_{(a,b]}d\rho_1^{(n)}$ to
$\int_{(a,b]}d\rho_1$ follows from the weak convergence of
$\rho_1^{(n)}$ to $\rho_1$. The pointwise convergence of $J^{(n)}$
to $J$ has been shown above. This implies the statement of the
lemma.
\end{proof}

Also in this framework the collapsing operators $\mathbb C_k$ are
defined likewise in \eqref{cat1} and \eqref{cat2} in the cases of
particles configurations. The continuity properties of $\mathbb C_k$
are derived also from lemma \ref{lemmino}.

\section{Empirical measures}

\label{collaempmis}
Given a configuration $\eta\in X_N$ we associate
to it its empirical measure $\pi_N(\eta)$. This is an element of
$\mathcal M$ defined as
$$
\pi_N(\eta):=\frac{1}{N}\sum_{x\in \mathbb
Z_N}\eta(x)\delta_{\frac{x}{N}}\ .
$$
Given a collection $(\eta_1,\dots ,\eta_k)$ of $k$ configurations of
$X_N$ we will write
$$
\pi_N^k(\eta_1,\dots ,\eta_k):=\left(\pi_N(\eta_1), \dots
,\pi_N(\eta_k)\right)\ .
$$
It is easy to check that
\begin{equation}
\pi_N\left(C_{\eta_2}[\eta_1]\right)=C_{\pi_N(\eta_2)}[\pi_N(\eta_1)]\
. \label{commuta}
\end{equation}
This is an important identity that justifies the fact that we used
the same symbols for the collapsing operators in different
frameworks. Note in fact that in \eqref{commuta} we are using the
same symbol $C$ with different meanings. Equation \eqref{commuta} is
the key identity to check the following commutation property
\begin{equation}
\pi_N^2\circ \mathbb C= \mathbb C\circ \pi_N^2\ .
\label{supercommuta}
\end{equation}
The validity of \eqref{supercommuta} directly implies its
generalization
\begin{equation}
\pi_N^k\circ \mathbb C_k= \mathbb C_k\circ \pi_N^k\ .
\label{supercommutak}
\end{equation}

Note that the same commutation relations hold also in the case we
had given a slightly different definition of empirical measure
associated to a configuration of particles of the TASEP. Sometimes
the empirical measure associated to a configuration $\eta\in X_N$ is
defined as an element of $\mathcal M^{0,b}$ whose density is a.e.
$$
\pi_N(\eta)(u)=\sum_{x\in \mathbb
Z_N}\eta(x)\chi_{_{\left[\frac{x}{N}-\frac{1}{2N},
\frac{x}{N}+\frac{1}{2N}\right)}}(u)\ .
$$
Also in this case \eqref{commuta}, \eqref{supercommuta} and
\eqref{supercommutak} hold.

In the case of the HAD process there is not a natural scale
parameter as in the case of the TASEP where it is the size of the
lattice. We will consider families, with index a natural number $N$,
of HAD models containing $[Nm]$ particles, with $m$ a positive real
number and $[\cdot ]$ the integer part. The scale parameter is $N$
and for any configuration $\underline x\in \Omega$ we define the
empirical measure as
$$
\pi_N(\underline x):=\frac{1}{N}\sum_i\delta_{x_i}\ .
$$
Note that we do not require in this definition that $N$ coincides
with $|\underline x|$. As in the case of the TASEP, given a
collection $(\underline{x}^{(1)},\dots ,\underline{x}^{(k)})$ of $k$
configurations of $\Omega$ we will write
$$
\pi_N^k(\underline{x}^{(1)},\dots
,\underline{x}^{(k)}):=\left(\pi_N(\underline{x}^{(1)}), \dots
,\pi_N(\underline{x}^{(k)})\right)\ .
$$
Also in this framework the analogous of \eqref{commuta},
\eqref{supercommuta} and \eqref{supercommutak} hold.

\section{Large deviations for uniform distributions}
\label{ldpu}

In this section we quickly derive, from well known results, large
deviations principles for the invariant measures of the TASEP and
the HAD process with particles of only one class.

\vskip0.5cm

Given a sequence $\mu_N$ of probability measures on a Polish metric
space $X$ we say that it satisfies a large deviation principle (LDP)
with parameter $N$ and rate function $I:X\to \mathbb R^+\cup
\left\{+\infty\right\}$, if
$$
\limsup_{N\to +\infty}N^{-1}\log \mu_N(C)\leq -\inf_{x\in C}I(x) \ \
\ \ \ \ \  \forall C\subseteq X \ , \ C \ closed \ ;
$$
$$
\liminf_{N\to +\infty}N^{-1}\log \mu_N(O)\geq -\inf_{x\in O}I(x) \ \
\ \ \ \ \ \forall O\subseteq X \ , \ O \ open \ .
$$
The rate function $I$ is lower semicontinuous and is called
\emph{good} if it has compact level sets.

\vskip0.5cm

Let us recall a classical LDP result for sampling without
replacement as stated for example in \cite{DZ}.

\noindent Let $\underline y^{(N)}=\left\{y_1^{(N)}, \dots,
y_N^{(N)}\right\}$ be $N$ elements of $\Lambda$ such that
$$
\pi_N(\underline
y^{(N)})=\frac{1}{N}\sum_{i=1}^{N}\delta_{y_i^{(N)}}\stackrel{N\to
+\infty}{\to}\lambda\ ,
$$
where $\to$ is the weak convergence and $\lambda \in \mathcal M$.
Let $\nu_N^M$ be the uniform measure on $M$-uples of elements of
$\underline y^{(N)}$, $M\leq N$:
\begin{equation}
\nu_N^M\left(\left\{y_{i_1}^{(N)}, \dots
,y_{i_M}^{(N)}\right\}\right)=\left\{
\begin{array}{cc}
\binom{N}{M}^{-1} & if \ \ i_{j}\neq i_{k} \ \ \forall j\neq k\ ,\\
0& otherwise\ .\\
\end{array}
\right. \nonumber
\end{equation}
This is the measure obtained from a procedure of $M$ samplings
without replacement among the $N$ elements of $\underline y^{(N)}$.
Consider now the case $M=M(N)=[Nm]$, where $[\ \cdot \ ]$ is the
integer part and $m\in (0,1]$. Then when $N$ diverges the measures
$\nu_N^{M(N)}\circ \pi_N^{-1}$ satisfy a LDP on $\mathcal M$
equipped with the weak topology with parameter $N$ and with a good
and convex rate function given by
\begin{equation}
\left\{
\begin{array}{cc}
mH\left(\frac{\rho}{m}\Big|\lambda\right)+(1-m)H\left(\frac{\lambda-\rho}{1-m}\Big|\lambda\right)
& if \ \
\lambda-\rho\in \mathcal M_{1-m}\ ,\\
+\infty & otherwise \ .
\end{array}
\right. \label{arrayone}
\end{equation}
Where $H\left( \cdot | \cdot  \right)$ is the relative entropy.

To get the rate functional for the invariant measure of the TASEP
with $[Nm]$ particles we need to consider $y_i^{(N)}=\frac{i}{N}$.
With this choice the measure $\lambda$ coincides with the Lebesgue
measure and consequently
\begin{equation}
H(\rho|\lambda)=\left\{
\begin{array}{ll}
\int_{\Lambda}\rho(u)\log \rho(u) du & if \ \rho\in \mathcal M^0\ ,\\
+\infty & otherwise\ .
\end{array}
\right.\nonumber
\end{equation}
In this case the condition for finiteness in \eqref{arrayone} is
equivalent to $\rho\in \mathcal M_m^{0,b}$. Observing that in this
case the collection of all the M-uples of $\underline y^{(N)}$ is in
bijection with $X_{N,M}$ (using the correspondence $\eta(x)=1$ if
and only if $\frac{x}{N}$ belongs to the M-uple), we obtain the
following result.

\begin{Pro}
\label{bernoulli} Let $\nu_N^M$ be the invariant measure of the
TASEP with $M$ particles. When N diverges the family of measures
$\nu_N^{[Nm]}\circ \pi_N^{-1}$ satisfy a LDP on $\mathcal M$
equipped with weak topology, with parameter $N$ and with good and
convex rate function
\begin{equation}
S_1(\rho):=\left\{
\begin{array}{cc}
\int_{\Lambda} du \ h_m(\rho(u)) & if \ \ \rho\in \mathcal M_m^{0,b}\ ,\\
+\infty & otherwise \ ,
\end{array}
\right. \nonumber
\end{equation}
where $h_m(x)=x\log\frac{x}{m}+(1-x)\log\frac{1-x}{1-m}$.
\end{Pro}
We use the symbol $S_1$ to remark that it is the rate function for a
one class system. The parameter $m$ is understood.

A LDP for the invariant measures of the HAD process with particles
of only one class follows directly from Sanov theorem (see \cite{DZ}
section 6.2). Sanov theorem states that if $X_i$ are i.i.d random
variables taking values in $\Lambda$ and having common law $\lambda$
then the empirical measure
$$
\frac{1}{N}\sum_{i}^{N}\delta_{X_i}
$$
satisfies a LDP on $\mathcal M$ equipped with the weak topology with
parameter $N$ and with a good convex rate functional given by the
relative entropy
\begin{equation}
\left\{
\begin{array}{lc}
H\left(\rho|\lambda\right) & if \ \rho\in \mathcal M_1\ ,\\
+\infty & otherwise\ .\\
\end{array}
\right. \nonumber
\end{equation}

From this general fact we can easily deduce the following
proposition as a special case when $\lambda$ is the Lebesgue
measure.
\begin{Pro}
\label{bernoulli2} Let $\mu_M$ be the invariant measure of the HAD
process with $M$ points and let $m$ be a positive real number. When
$N$ diverges the family of measures $\mu_{[Nm]}$ satisfy a large
deviation principle on $\mathcal M$ equipped with the weak topology
with parameter $N$ and with the good and convex rate functional
\begin{equation}
S_1(\rho):=\left\{
\begin{array}{cc}
\int_{\Lambda} du \ k_m(\rho(u)) & if \ \ \rho\in \mathcal M^0_m\ ,\\
+\infty & otherwise\ ,
\end{array}
\right. \nonumber
\end{equation}
where $k_m(x)=x\log\frac{x}{m}$.
\end{Pro}
For simplicity we use the same symbol $S_1$ already used for the
TASEP, but the rate functions are different.

\section{LDP for 2-class processes}
\label{ldp2}

Theorem \ref{bernoulli} is immediately generalized to the case of
product measures. Consider $0<m_1\leq m_2\leq 1$ real numbers. The
family of measures
\begin{equation}
\left(\nu_N^{[Nm_1]}\times \nu_N^{[Nm_2]}\right)\circ
\left(\pi_N^{2}\right)^{-1}\label{met2}
\end{equation}
satisfy a LDP on $\mathcal M \times \mathcal M$ endowed with weak
topology, with parameter $N$ and with good and convex rate function
given by
\begin{equation}
\left\{
\begin{array}{cc}
\sum_{i=1}^2\int_{\Lambda} du \ h_{m_i}(\rho_i(u)) & if \ \ \rho_i\in \mathcal M^{0,b}_{m_i}\ ,\\
+\infty & otherwise\ .
\end{array}
\right. \label{met2func}
\end{equation}
We are interested in proving a LDP for the empirical measures of the
invariant measures of the 2-class TASEP. This means that we are
interested in proving a LDP for the sequence of measures
\begin{equation}
\left[\left(\nu_N^{[Nm_1]}\times \nu_N^{[Nm_2]}\right)\circ
\left(\mathbb C_{2}\right)^{-1}\right]\circ
\left(\pi_N^2\right)^{-1}\ ,
%\label{chissa}
\end{equation}
that due to identity \eqref{supercommuta} coincides with the
sequence of measures
\begin{equation}
\left[\left(\nu_N^{[Nm_1]}\times \nu_N^{[Nm_2]}\right)\circ
\left(\pi_N^{2}\right)^{-1}\right]\circ \left(\mathbb
C_2\right)^{-1}\ . \label{chissa2}
\end{equation}
Lemma \ref{lemmino} suggests that we can apply the contraction
principle. \vskip0.5cm

The contraction principle (see for example \cite{DZ} section 4.2.1)
states that if $\mu_N$ is a sequence of measures satisfying a large
deviation principle on a Polish metric space $X$ with a good rate
functional $I(x)$ and $f:X\to Y$ is a continuous map from X to
another Polish metric space $Y$, then also the sequence of measures
$\mu_N\circ f^{-1}$ satisfy a LDP on $Y$ with good rate functional
$K$ given by $K(y)=\inf_{\left\{x\in f^{-1}(y)\right\}}I(x)$. This
formulation can in fact be extended (see remark c at page 127 of
\cite{DZ}) to the case when the map $f$ is continuous only at the
elements $x\in X$ such that $I(x)<+\infty$. This is exactly our
setting.

\vskip0.5cm

Before stating the theorem that we obtain following this strategy we
need some facts and notations. Consider an interval $\mathcal
U:=[u^l,u^r]\subset \Lambda$ and a density $\rho(u)$. Let us
introduce the extended function
\begin{equation}
F_{\rho}^{\mathcal U}(u):=\left\{
\begin{array}{cc}
\int_{u^l}^u\rho (w) dw & if \ \ u\in \mathcal U\ ,\\
-\infty & otherwise\ .\\
\end{array}
\right.\nonumber
\end{equation}
Note that in the case of measures in $\mathcal M^0$ we can use the
above notation for integration because there is no difference
between open and closed intervals. Let us also consider the extended
function $\mathbb F_{\rho}^{\mathcal U}$ on the real line $\mathbb
R$ defined as
\begin{equation}
\mathbb F_{\rho}^{\mathcal U}(u):=\left\{
\begin{array}{cc}
F_{\rho}^{\mathcal U}(u^l+u) & if \ \ u\in [0,|[u^l,u^r]|]\ ,\\
-\infty & otherwise\ .\\
\end{array}
\right.\nonumber
\end{equation}
We call $\widehat{\mathbb F}_{\rho}^{\mathcal U}$ the concave
envelope of $\mathbb F_{\rho}^{\mathcal U}(u)$. We then define the
following extended function on $\Lambda$
\begin{equation}
\widehat{F}_{\rho}^{\mathcal U}(u):= \widehat{\mathbb
F}_{\rho}^{\mathcal U}\left(|[u^l,u]|\right) \nonumber
\end{equation}

Note that the following properties hold
\begin{equation}
\left\{
\begin{array}{l}
\widehat{F}_{\rho}^{\mathcal U}(u)\geq F_{\rho}^{\mathcal U}(u)\ ,\\
\widehat{F}_{\rho}^{\mathcal U}(u^l)=F_{\rho}^{\mathcal U}(u^l)=0\ ,\\
\widehat{F}_{\rho}^{\mathcal U}(u^r)=F_{\rho}^{\mathcal U}(u^r)=\int_{u^l}^{u^r}\rho(w)dw\ .\\
\end{array}
\right.\nonumber
\end{equation}
Finally we call $\rho^{\mathcal U}$ a density on $\Lambda$ such that
$\int_{u^l}^{u} \rho^{\mathcal U}(w)dw=\widehat{F}_{\rho}^{\mathcal
U}(u)$ for any $u\in \mathcal U$ and $\rho^{\mathcal
U}(u)\chi_{_{\mathcal U^c}}(u)=0$ a.e.. It can be shown that
$\rho^{\mathcal U}$ is a positive measure and it belongs to
$\mathcal M^{0,b}$ when $\rho \in \mathcal M^{0,b}$.

Given $\rho_1$ and $\rho_2$ measurable functions on $\Lambda$, the
set $\mathcal U:=\left\{u\in \Lambda : \
\rho_1(u)=\rho_2(u)\right\}$ is a.e. equivalent to the disjoint
union of at most countable many closed intervals $\mathcal
U_i:=\left[u_i^l,u_i^r\right]$. Indeed it is simple to check that
the family of sets with this property is a $\sigma$-algebra that
includes all the open sets. Hence it must include all the Borel
sets.

\begin{Le}
\label{gigolino}
Let $(\rho_1,\rho_2)\in I^{2,\uparrow}\cap
\left(\mathcal M_{m_1}^{0}\times \mathcal M_{m_2}^{0}\right)$. We
have that the pair $(\psi_1,\psi_2)\in \mathcal M_{m_1}^{0}\times
\mathcal M_{m_2}^{0}$ is such that $\mathbb
C(\psi_1,\psi_2)=(\rho_1,\rho_2)$ if and only if the following
conditions are satisfied
\begin{equation}
\left\{
\begin{array}{l}
\psi_2=\rho_2 \ ,\\
\psi_1(u)\chi_{_{\mathcal U^c}}(u)=\rho_1(u)\chi_{_{\mathcal U^c}}(u)\ , \ \ a.e. \ , \\
F_{\psi_1}^{\mathcal U_i}(u^r_i)= F_{\rho_2}^{\mathcal U_i}(u^r_i) \
, \ \ \
\forall i \ ,\\
F_{\psi_1}^{\mathcal U_i}\geq F_{\rho_2}^{\mathcal U_i} \ , \ \ \
\forall
i\ .\\
\end{array}
\right. \label{stracond}
\end{equation}
\end{Le}
\begin{proof}
We recall that due to the fact that the measures are absolutely
continuous with respect to Lebesgue measure, the above functions $F$
are continuous on the interior part of the intervals where they are
different from $-\infty$. First we show that given a pair
$(\psi_1,\psi_2)\in \mathcal M_{m_1}^{0}\times \mathcal M_{m_2}^{0}$
that satisfies conditions \eqref{stracond} then we have $\mathbb
C(\psi_1,\psi_2)=(\rho_1,\rho_2)$. Clearly we need only to prove
that $C_{\rho_2}[\psi_1]=\rho_1$. Remember that for measures
belonging to $\mathcal M^0$ the collapsing procedure acts as in
\eqref{azasscont}. Conditions number three and four in
\eqref{stracond} imply that for any $u\in \mathcal U_i$ and for any
$i$ we have
\begin{equation}
\int_{u}^{u_i^r}d\psi_1- \int_{u}^{u_i^r}d\rho_2\leq 0\ . \nonumber
\end{equation}
Let us consider the following subset of $\mathcal U$
\begin{equation}
\widetilde{\mathcal U}:=\bigcup_i\left\{u\in \mathcal U_i : \
F_{\psi_1}^{\mathcal U_i}(u)> F_{\rho_2}^{\mathcal U_i}(u) \right\}\
. \nonumber
\end{equation}
From the fact that equality between the $F$ holds at the boundary of
each $\mathcal U_i$ and from continuity in the interior part we
deduce that this is an open set. We show now that
$\widetilde{\mathcal U}$ coincides with the set $\mathcal J$ as
defined in equation \eqref{dafare} for the pair of measures
$(\psi_1,\rho_2)$. Clearly $\widetilde{\mathcal U}\subseteq \mathcal
J$. This follows from the fact that for any $u\in
\widetilde{\mathcal U}$ with $u\in \mathcal U_{i}$ it holds
$$
\int^{u}_{u_i^l}d\psi_1- \int^{u}_{u_i^l}d\rho_2> 0\ .
$$
To prove that $\widetilde{\mathcal U}=\mathcal J$ we need to show
that for any $u\in \widetilde{\mathcal U}^c$ and for any $v\in
\Lambda$ it holds
$$
\int^{u}_{v}d\psi_1- \int^{u}_{v}d\rho_2\leq 0\ .
$$
Note that
$$
\int_{[v,u]\cap \mathcal U^c}d\psi_1- \int_{[v,u]\cap \mathcal
U^c}d\rho_2\leq 0\ ,
$$
due to the fact that on $\mathcal U^c$ we have a.e. that
$\psi_1(w)=\rho_1(w)$ and $\rho_1(w)<\rho_2(w)$. We consider first
the case $u\in\mathcal U^c$. In this case we have
\begin{eqnarray}
& & \int^{u}_{v}d\psi_1- \int^{u}_{v}d\rho_2 = \int_{[v,u]\cap
\mathcal U^c}d\psi_1- \int_{[v,u]\cap \mathcal U^c}d\rho_2\nonumber
\\
&+ & \sum_{\mathcal U_i\subseteq [v,u]}\left(\int_{\mathcal
U_i}d\psi_1- \int_{\mathcal U_i}d\rho_2\right)
+\left(\int_{v}^{u_{i^*}^r}d\psi_1-
\int_{v}^{u_{i^*}^r}d\rho_2\right)\chi_{_{\mathcal U}}(v)\ .
\label{ehehvit}
\end{eqnarray}
When $v\in \mathcal U$ we called $\mathcal U_{i^*}$ the interval to
which it belongs. All terms on the right hand side of
\eqref{ehehvit} are nonpositive. Consider now the case $u\in
\mathcal U \cap \widetilde{\mathcal U}^c$ and call $\mathcal
U_{j^*}$ the interval to which it belongs. In this case we need to
modify formula \eqref{ehehvit} multiplying the last term on the
right hand side by $\chi_{_{[v,u]}}(u_{j^*}^l)$ and adding
$$
\left(\int_{v}^{u}d\psi_1-
\int_{v}^{u}d\rho_2\right)\chi_{_{[v,u]^c}}(u_{j^*}^l)+
\left(\int^{u}_{u_{j^*}^l}d\psi_1-
\int^{u}_{u_{j^*}^l}d\rho_2\right)\chi_{_{[v,u]}}(u_{j^*}^l)\ .
$$
Both terms are nonpositive. Note that $i^*$ can coincide with $j^*$.

Conversely we assume that for $(\psi_1,\psi_2)\in \mathcal
M_{m_1}^{0}\times \mathcal M_{m_2}^{0}$ holds $\mathbb
C(\psi_1,\psi_2)= (\rho_1,\rho_2)$ and show that this implies the
validity of conditions \eqref{stracond}.

The validity of the first condition in \eqref{stracond} is obvious
due to the fact that the collapsing operator $\mathbb C$ preserves
the second component.

Using \eqref{azasscont} and the previous statement we have that a.e.
it holds the following equality
$$
C_{\rho_2}[\psi_1](u)=\psi_1(u)\chi_{_{\mathcal J^c}}(u)+
\rho_2(u)\chi_{_{\mathcal J}}(u)=\rho_1(u)\ .
$$
We multiply both sides by $\chi_{_{\mathcal U^c}}(u)$ and obtain
\begin{equation}
\left(\psi_1(u)-\rho_1(u)\right)\chi_{_{\mathcal J^c\cap \mathcal
U^c}}(u)+ \left(\rho_2(u)-\rho_1(u)\right)\chi_{_{\mathcal J\cap
\mathcal U^c}}(u)=0\ , \ \ \ \ a.e.\ .\label{piangesempre}
\end{equation}
Due to the fact that the two terms can be different from zero on
disjoint sets we have that \eqref{piangesempre} is equivalent to the
two equations
\begin{equation}
\left(\psi_1(u)-\rho_1(u)\right)\chi_{_{\mathcal J^c\cap \mathcal
U^c}}(u)=0\ , \ \ \ \ a.e.\ ,\label{piangesempre1}
\end{equation}
\begin{equation}
\left(\rho_2(u)-\rho_1(u)\right)\chi_{_{\mathcal J\cap \mathcal
U^c}}(u)=0\ , \ \ \ \ a.e.\ .\label{piangesempre2}
\end{equation}
From the fact that $\rho_1(u)<\rho_2(u)$ a.e. on $\mathcal U^c$ and
from equation \eqref{piangesempre2} we deduce that
\begin{equation}
\int_{\mathcal J^c\cap \mathcal U^c}d\lambda=\int_{\mathcal
U^c}d\lambda\ , \label{infinitilabel}
\end{equation}
where $\lambda$ is the Lebesgue measure. Equation
\eqref{piangesempre1} imposes that $\psi_1(u)=\rho_1(u)$ a.e. on
$\mathcal J^c\cap \mathcal U^c$ and using \eqref{infinitilabel} we
obtain the validity of the second condition in \eqref{stracond}.

In the case of absolutely continuous measures we have that $\mathcal
J$ is an open set and this implies that $u_i^l$ and $u_i^r$ do not
belong to $\mathcal J$ for any $i$. In fact let us suppose for
example that $u_i^r\in \mathcal J$, then there exists a ball
$B_{u_i^r}(\epsilon)$ centered in $u_i^r$ such that
$B_{u_i^r}(\epsilon)\subseteq \mathcal J$ and
$\int_{B_{u_i^r}(\epsilon)\cap \mathcal U^c}d\lambda >0$ (recall
that $\lambda$ denotes Lebesgue measure). This is impossible due to
\eqref{infinitilabel}. From the fact that $u_i^l$ and $u_i^r$ belong
to $\mathcal J^c$ we deduce that
$$
J(u_i^l)=J(u_i^r)=0\ .
$$
Using this we obtain
\begin{equation}
\int_{u_i^l}^{u_i^r}d\rho_2=\int_{u_i^l}^{u_i^r}dC_{\rho_2}[\psi_1]=\int_{u_i^l}^{u_i^r}d\psi_1
\ , \label{cricri}
\end{equation}
that is the third condition in \eqref{stracond}.

Let us suppose that there exists an $u\in \mathcal U_i$ such that
$F_{\psi_1}^{\mathcal U_i}(u)<F_{\rho_2}^{\mathcal U_i}(u)$. Then
using also \eqref{cricri} we obtain
$$
\int_{u}^{u_i^r}d\psi_1-\int_{u}^{u_i^r}d\rho_2=F_{\rho_2}^{\mathcal
U_i}(u)-F_{\psi_1}^{\mathcal U_i}(u)>0\ .
$$
This implies $u_i^r\in \mathcal J$ that is impossible. As a
consequence we obtain the validity of the fourth condition in
\eqref{stracond}.
\end{proof}

Now we can state and prove our large deviations result.

\begin{Th}
\label{Th1} Let $0<m_1< m_2\leq 1$ real numbers. Consider the
$2$-class TASEP on $\mathbb Z_N$ with respectively $[Nm_1]$ first
class particles and $[Nm_2]$ total particles. When the pair
$(\eta_1,\eta_2)$ is distributed according to the invariant measure
of the process, we have that $\pi_N^2(\eta_1,\eta_2)$ satisfies a
LDP with parameter $N$ and with good rate function
$S_2(\rho_1,\rho_2)$ defined as follows. It takes the value
$+\infty$ when $(\rho_1,\rho_2)\not \in
I^{2,\uparrow}\cap\left(\mathcal M^{0,b}_{m_1}\times \mathcal
M^{0,b}_{m_2}\right)$. When $(\rho_1,\rho_2) \in
I^{2,\uparrow}\cap\left(\mathcal M^{0,b}_{m_1}\times \mathcal
M^{0,b}_{m_2}\right)$ it takes the value
\begin{equation}
S_2(\rho_1,\rho_2)= \int_{\mathcal U^c}h_{m_1}(\rho_1(u))du
+\sum_i\int_{\mathcal U_i}h_{m_1}\left(\rho_1^{\mathcal
U_i}(u)\right)+\int_{\Lambda}h_{m_2}(\rho_2(u))du\ ,
\end{equation}
where $ \mathcal U_i$ are disjoint closed intervals and $\mathcal
U=\cup_i \mathcal U_i$ is a.e. equivalent to the set
$$
\left\{u\in \Lambda : \ \rho_1(u)=\rho_2(u) \right\}\ .
$$
The symbol $\rho_1^{\mathcal U_i}$ is defined before the statement
of lemma \ref{gigolino}.
\end{Th}
\begin{proof}
As outlined before we can apply the contraction principle. The
sequence of measures in \eqref{met2} satisfies a LDP with good rate
functional given by \eqref{met2func} and, as shown in lemma
\ref{lemmino}, the map $\mathbb C$ is continuous on every point
where \eqref{met2func} is different from $+\infty$. We immediately
get that the sequence of measures \eqref{chissa2} satisfies a LDP
with a good rate given by
\begin{equation}
S_2(\rho_1,\rho_2)=\inf_{\left\{(\psi_1,\psi_2)\in \mathcal
M^{0,b}_{m_1}\times M^{0,b}_{m_2} : \ \mathbb
C(\psi_1,\psi_2)=(\rho_1,\rho_2)\right\}}\left(\int_{\Lambda}h_{m_1}(\psi_1(u))du+
\int_{\Lambda}h_{m_2}(\psi_2(u))du\right)\ . \label{servira}
\end{equation}
By convention the infimum over an empty set is defined as $+\infty$
and remember that we are calling $\psi_i$ both the measures and the
corresponding densities. From \eqref{servira} we see immediately
that $S_2$ is equal to $+\infty$ on
$\left(I^{2,\uparrow}\cap\left(\mathcal M^{0,b}_{m_1}\times \mathcal
M^{0,b}_{m_2}\right)\right)^c$.

\vskip0.5cm

The case $m_1=m_2=m$ is not covered by the theorem (in fact it could
be, choosing appropriately an interval $\mathcal U$). In this case
the rate functional is different from $+\infty$ only when
$\rho_1=\rho_2$ and consequently $\mathcal U=\Lambda$. The rate
functional is then
\begin{equation}
S_2(\rho_1,\rho_2)=\left\{
\begin{array}{lc}
\int_{\Lambda}h_{m}(\rho(u))du & if \ \rho_1=\rho_2=\rho \in
\mathcal M_m^{0,b}\ ,\\
+\infty & otherwise\ .\\
\end{array}
\right. \nonumber
\end{equation}
When $m_1<m_2$ we have a strict inclusion $\mathcal U\subset
\Lambda$.

\vskip0.5cm

First we prove existence of a minimizer for the variational problem
\eqref{servira}, then we prove uniqueness and finally we
characterize it.

In the case of the TASEP, existence of a minimizer can be derived
directly showing that \eqref{servira} is a minimization problem of a
lower semicontinuous functional over a compact set.

We show instead a more involved proof that works also in the case of
the HAD process. When $(\rho_1,\rho_2)\in I^{2,\uparrow}$ then
$\mathbb C(\rho_1,\rho_2)=(\rho_1,\rho_2)$. This means that the we
can modify the infimum in \eqref{servira} restricting to pairs
$(\psi_1,\psi_2)$ that satisfy the further condition
\begin{equation}
\left\{(\psi_1,\psi_2)\in \mathcal M^{0,b}_{m_1}\times M^{0,b}_{m_2}
: \  \sum_{i=1}^2\int_{\Lambda}h_{m_i}(\psi_i(u))du \leq
\sum_{i=1}^2\int_{\Lambda}h_{m_i}(\rho_i(u))du \right\}\ .
\label{compatto}
\end{equation}
Due to the fact that \eqref{met2func} is a good rate function we
have that \eqref{compatto} is a compact set. The constraints set in
\eqref{servira} is easily seen to be a closed subset of $\mathcal
M_{m_1}^{0,b}\times \mathcal M_{m_2}^{0,b}$ (and this holds also
without the condition of bounded densities). In fact consider a
sequence $(\psi_1^{(n)},\psi_2^{(n)})$ belonging to this set and
converging to $(\psi_1,\psi_2)\in \mathcal M_{m_1}^{0,b}\times
\mathcal M_{m_1}^{0,b}$. Then by lemma \ref{lemmino} we have that
$$
\mathbb C(\psi_1,\psi_2)= \mathbb C\left(\lim_{n\to
\infty}(\psi_1^{(n)},\psi_2^{(n)})\right)=\lim_{n\to \infty}\mathbb
C(\psi_1^{(n)},\psi_2^{(n)})=(\rho_1,\rho_2)\ ,
$$
which means that the set is closed. We obtained an infimum of a
lower semicontinuous functional over a compact set and the existence
of a minimizer follows.

We prove now uniqueness of the minimizer. The functional to be
minimized in \eqref{servira} is strictly convex. This follows
directly from the fact that the real functions $h_{m_i}$ are
strictly convex. The set on which we are minimizing is also a convex
subset of $\mathcal M_{m_1}^{0,b}\times \mathcal M_{m_2}^{0,b}$.
This follows directly from its characterization given in lemma
\ref{gigolino}. If $(\psi_1,\psi_2)$ and $(\psi_1^*,\psi_2^*)$
satisfy conditions \eqref{stracond} then clearly also the convex
combination
$(\psi_{1c},\psi_{2c}):=c(\psi_1,\psi_2)+(1-c)(\psi_1^*,\psi_2^*)$
satisfies the same conditions and this is clearly true also for the
additional condition of bounded density. A classical result in
convex analysis \cite{ET} guarantees uniqueness of the minimizer of
a strictly convex functional over a convex set.

Finally we characterize the unique minimizer.

Using conditions \eqref{stracond} we write \eqref{servira} as
\begin{eqnarray*}
& S_2(\rho_1,\rho_2)=\int_{\Lambda}h_{m_2}(\rho_2(u))du+
\int_{\mathcal U^c}h_{m_1}(\rho_1(u))du\\
& +\sum_i\inf_{\left\{\phi_i\in  A_i\right\}}\int_{\mathcal
U_i}h_{m_1}(\phi_i(u))du\ ,
\end{eqnarray*}
where
\begin{equation}
A_i:={\left\{\phi\in \mathcal M^{0,b}: \ F_{\phi}^{\mathcal U_i}\geq
F_{\rho_2}^{\mathcal U_i}; \ F_{\phi}^{\mathcal U_i}(u^r_i)=
F_{\rho_2}^{\mathcal U_i}(u^r_i)\right\}}\ . \label{confuso}
\end{equation}
We need then to study the variational problems
\begin{equation}
\inf_{\left\{\phi_i\in A_i\right\}}\int_{\mathcal
U_i}h_{m_1}(\phi_i(u))du\ . \label{labella}
\end{equation}
Existence and uniqueness of the minimizers can be shown as before.
It remains to characterize them.

From the strict convexity in $x$ of the real function $h_m(x)$ and
Jensen inequality we have for any interval $[v_1,v_2]\subseteq
\Lambda$ and any density $\phi$
\begin{equation}
\frac{1}{|[v_1,v_2]|}\int_{v_1}^{v_2}h_m(\phi(u))du\geq
h_m\left(\frac{1}{|[v_1,v_2]|}\int_{v_1}^{v_2}\phi(u)du\right)\ .
\label{geojensen}
\end{equation}
Moreover this inequality is strict as soon as the density $\phi(u)$
is not a.e. constant. A geometric interpretation of this inequality
is the following. Let $\phi$ and $\phi'$ be two densities defined on
an interval $\mathcal V\subset \Lambda$ and let $[v_1,v_2]\subseteq
\mathcal V$. Consider the case in which $F_{\phi}^{\mathcal
V}(u)=F_{\phi'}^{\mathcal V}(u)$ for any $u\not \in (v_1,v_2)$ and
the graph of $F_{\phi'}^{\mathcal V}$ when $u\in [v_1,v_2]$ linearly
interpolates $(v_1,F_{\phi}^{\mathcal V}(v_1))$ and
$(v_2,F_{\phi}^{\mathcal V}(v_2))$. More precisely
\begin{equation}
F_{\phi'}^{\mathcal V}(u)=\frac{F_{\phi}^{\mathcal
V}(v_2)-F_{\phi}^{\mathcal
V}(v_1)}{|[v_1,v_2]|}(|[v_1,u]|)+F_{\phi}^{\mathcal V}(v_1)\ ; \ \ \
\ \ \ u\in [v_1,v_2]\ . \label{picco}
\end{equation}
Note that if $\phi\in \mathcal M^{0,b}$ then necessarily also
$\phi'\in \mathcal M^{0,b}$. Inequality \eqref{geojensen} then
simply says that
\begin{equation}
\int_{\mathcal V}h_m(\phi(u))du \geq \int_{\mathcal
V}h_m(\phi'(u))du\ , \label{muffa}
\end{equation}
with the strict inequality holding if $\phi(u)$ and $\phi'(u)$ do
not coincide a.e.. In the rest of the proof we consider pairs
$(\phi,\phi')$ whose corresponding $F_{\phi}^{\mathcal U_i}(u)$ and
$F_{\phi'}^{\mathcal U_i}(u)$ are related as before. We can then
apply inequality \eqref{muffa}.

Consider a $\phi\in  A_i$, such that $\mathbb F_{\phi}^{\mathcal
U_i}$ is not a concave function. Then clearly there exists a $\phi'$
such that $F_{\phi'}^{\mathcal U_i}(u)\geq F_{\phi}^{\mathcal
U_i}(u)$ for any $u$. In particular this implies that also $\phi'\in
A_i$. Inequality \eqref{muffa} implies that $\phi$ can not be the
unique minimizer. The unique minimizer $\phi_i^*$ of \eqref{labella}
has then necessarily $\mathbb F_{\phi_i^*}^{\mathcal U_i}$ concave.

Let us now consider $\phi, \psi \in A_i$ such that both $\mathbb
F_{\phi}^{\mathcal U_i}$ and $\mathbb F_{\psi}^{\mathcal U_i}$ are
concave and moreover there exists an $u\in \mathcal U_i$ such that
$\mathbb F_{\phi}^{\mathcal U_i}(u)>\mathbb F_{\psi}^{\mathcal
U_i}(u)$. Then $\phi$ can not be the unique minimizer. We can in
fact construct a $\phi'$ considering an affine function through
$(u,\mathbb F_{\psi}^{\mathcal U_i}(u))$, obtained using an element
in the superdifferential of $\mathbb F_{\psi}^{\mathcal U_i}$ at
$u$. We then have
$$
F_{\phi'}^{\mathcal U_i}(v)\geq \min\left\{F_{\phi}^{\mathcal
U_i}(v),F_{\psi}^{\mathcal U_i}(v)\right\}\ , \ \ \ \forall v\in
\mathcal U_i \ ,
$$
so that $\phi'\in A_i$. Inequality \eqref{muffa} then implies that
$\phi$ can not be the unique minimizer.
%$$
%\int_{\mathcal U_i}h_m(\phi(v))dv \geq \int_{\mathcal
%U_i}h_m(\phi'(v))dv\ .
%$$
The unique minimizer $\phi_i^*$ is then necessarily such that
$\mathbb F_{\phi^*_i}^{\mathcal U_i}$ is the smallest among all the
concave functions that are above $\mathbb F_{\rho_2}^{\mathcal
U_i}$, that is its \emph{concave envelope} $\widehat{\mathbb
F}_{\rho_2}^{\mathcal U_i}$. This shows that on the interval
$\mathcal U_i$ we have $\phi_i^*(u)=\rho_2^{\mathcal
U_i}(u)=\rho_1^{\mathcal U_i}(u) \ a.e.$.
\end{proof}

The rate functional $S_2$ is non negative and is zero if and only if
$\psi_1(u)=m_1 \ a.e.$ and $\psi_2(u)=m_2 \ a.e.$ and corresponding
$\rho_2(u)=m_2 \ a.e.$ and $\rho_1(u)=C_{m_2}[m_1](u)=m_1 \ a.e.$.

The rate functional $S_2$ is not convex. This can be shown by the
following example. Let us consider densities $(\rho_1,\rho_2)$ and
$(\rho_1^*,\rho_2^*)$ defined a.e. as
$$
\rho_1(u)=\chi_{_{[\frac{1}{4},\frac{1}{2}]}}(u)\ ; \ \ \
\rho_1^*(u)=\frac{1}{2}\chi_{_{[\frac{1}{2},1]}}(u)\ ;\ \ \
\rho_2(u)=\rho_2^*(u)=\chi_{_{[\frac{1}{4},1]}}(u)\ .
$$
We consider the convex combination
$$
(\rho_{1c},\rho_{2c}):=c(\rho_1,\rho_2)+(1-c)(\rho_1^*,\rho_2^*)\ .
$$
We have that
$$
S_2(\rho_1,\rho_2)=\int_{\Lambda}
h_{\frac{3}{4}}\left(\rho_2(u)\right) du
+\int_{[0,\frac{1}{2}]}h_{\frac{1}{4}}\left(\rho_1^{[0,\frac{1}{2}]}(u)\right)
du+ \int_{[\frac{1}{2},1]}h_{\frac{1}{4}}\left(\rho_1(u)\right) du\
,
$$
where
$\rho_1^{[0,\frac{1}{2}]}(u)=\frac{1}{2}\chi_{_{[0,\frac{1}{2}]}}(u)$
a.e.. We have also
$$
S_2(\rho_1^*,\rho_2^*)=\int_{\Lambda}
h_{\frac{3}{4}}\left(\rho_2(u)\right) du
+\int_{\Lambda}h_{\frac{1}{4}}\left(\rho_1^*(u)\right) du
$$
and for any $c\in(0,1)$
$$
S_2(\rho_{1c},\rho_{2c})=\int_{\Lambda}
h_{\frac{3}{4}}\left(\rho_2(u)\right) du
+\int_{\Lambda}h_{\frac{1}{4}}\left(\rho_{1c}(u)\right) du\ .
$$
Convexity of $S_2$ would imply the validity for any $c\in[0,1]$ of
the following inequality
\begin{equation}
cS_2(\rho_1,\rho_2)+(1-c)S_2(\rho_1^*,\rho_2^*)-S_2(\rho_{1c},\rho_{2c})\geq
0\ . \label{coxin}
\end{equation}
If we take the limit $c\uparrow 1$, on the left hand side of
\eqref{coxin} we obtain
$$
\int_{[0,\frac{1}{2}]}h_{\frac{1}{4}}(\rho_1^{[0,\frac{1}{2}]}(u))
du-\int_{[0,\frac{1}{2}]}h_{\frac{1}{4}}(\rho_1(u)) du\ ,
$$
that is clearly strictly negative. This shows that $S_2$ is not
convex.

\vskip0.5cm

To simplify notations in the following lemmas the fact that all the
measures are absolutely continuous, have bounded densities and have
a fixed total mass will be understood. It is understood also the
fact that $(\rho_1,\rho_2)\in I^{2,\uparrow}$.

We can still obtain interesting results from the contraction
principle. The following identity has to be satisfied
\begin{equation}
\inf_{\rho_1}S_2(\rho_1,\rho_2)=S_1(\rho_2)=\int_{\Lambda}h_{m_2}(\rho_2(u))du\
. \label{seilento}
\end{equation}
From the microscopic point of view this identity simply derives from
the fact that if we forget the labels first and second class and we
just look at positions of particles the dynamics that we observe is
a TASEP. From the variational point of view we have the following
lemma.
\begin{Le}
The unique minimizer $\rho_1^*$ in \eqref{seilento} such that
\begin{equation}
S_2(\rho_1^*,\rho_2)=S_1(\rho_2) \nonumber
\end{equation}
is given by
\begin{equation}
\rho_1^*=C_{\rho_2}[m_1]\ , \label{fff}
\end{equation}
where $m_1$ denotes the measure with a constant density equal to
$m_1$.
\end{Le}
\begin{proof}
We have
\begin{eqnarray*}
& \inf_{\left\{\rho_1\right\}}\inf_{\left\{(\psi_1,\psi_2): \
\mathbb
C(\psi_1,\psi_2)=(\rho_1,\rho_2)\right\}}\left(\int_{\Lambda}h_{m_1}(\psi_1(u))du+
\int_{\Lambda}h_{m_2}(\psi_2(u))du\right)\\
&= \int_{\Lambda}h_{m_2}(\rho_2(u))du+
\inf_{\psi_1}\int_{\Lambda}h_{m_1}(\psi_1(u))du\\
&=\int_{\Lambda}h_{m_2}(\rho_2(u))du\ .
\end{eqnarray*}
The unique minimizer in the expression above has a density
$\psi_1^*(u)=m_1 \ a.e.$ and correspondingly the unique minimizer
$\rho_1^*$ in \eqref{seilento} is obtained from the collapsing
procedure applied to $\psi_1^*$, that is \eqref{fff}.
\end{proof}

The minimizer $\rho_1^*$ in \eqref{fff} corresponds to the typical
density of first class particles when the system is conditioned to
have a total density $\rho_2$.

\vskip0.5cm

Still from the contraction principle we have that the following
identity has to be satisfied
\begin{equation}
\inf_{\rho_2}S_2(\rho_1,\rho_2)=S_1(\rho_1)=\int_{\Lambda}h_{m_1}(\rho_1(u))du\
. \label{pierpo}
\end{equation}
From the microscopic point of view this identity simply derives from
the fact that if we forget second class particles and observes only
first class particles what we see is a TASEP. Identity
\eqref{pierpo} can be deduced also from purely variational
arguments.

Fix the density $\rho_1$ and call $\mathcal V=\cup_i\mathcal V_i$ a
subset of $\Lambda$ a.e. equivalent to the subset $\left\{u\in
\Lambda: \ \rho_1(u)>m_2\right\}$. The sets $\mathcal V_i$ are
disjoint closed intervals, $\mathcal V_i:=[v_i^l,v_i^r]$. To any
such an interval we associate an element $w_i^l\in \Lambda$. This
element is the nearest $w\in \Lambda$ to the left of $v_i^l$ such
that
\begin{equation}
\int_w^{v_i^r}\left(m_2-\rho_1(z)\right) dz=0\ . \label{vaccino}
\end{equation}
It exists due to the fact that
$\int_{\Lambda}\left(m_2-\rho_1(z)\right) dz=m_2-m_1>0$, the
function of $w$ given by $\int_w^{v_i^r}\left(m_2-\rho_1(z)\right)
dz$ is continuous and condition \eqref{vaccino} identify a closed
set. Consider now the intervals $\widetilde{\mathcal
V}_i:=[w_i^l,v_i^r]$.

\begin{Le}
The unique minimizer $\rho_2^*$ in \eqref{pierpo} such that
\begin{equation}
S_2(\rho_1,\rho_2^*)=S_1(\rho_1)\ , \nonumber
\end{equation}
has a density a.e. equal to
\begin{equation}
\rho_2^*(u)=\left\{
\begin{array}{cc}
m_2 & if \ u\not \in \cup_i \widetilde{\mathcal V}_i\ ,\\
\rho_1(u) & if \ u\in \cup_i \widetilde{\mathcal V}_i\ .\\
\end{array}
\right.\label{ggg}
\end{equation}
\end{Le}
\begin{proof}
We have
\begin{eqnarray*}
& S_2(\rho_1,\rho_2)-S_1(\rho_1)=\int_{\Lambda}h_{m_2}(\rho_2(u))du+
\int_{\mathcal U^c}h_{m_1}(\rho_1(u))du\\
& +\sum_i\int_{\mathcal U_i}h_{m_1}\left(\rho_1^{\mathcal
U_i}(u)\right)-
\int_{\Lambda}h_{m_1}\left(\rho_1(u)\right)\\
&=\int_{\mathcal U^c}h_{m_2}(\rho_2(u))du+
\sum_i\int_{\mathcal U_i}h_{m_2}\left(\rho_1(u)\right)\\
& +\sum_i\int_{\mathcal U_i}h_{m_1}\left(\rho_1^{\mathcal
U_i}(u)\right)- \sum_i\int_{\mathcal
U_i}h_{m_1}\left(\rho_1(u)\right)\ .
\end{eqnarray*}
Where we used the fact that $\rho_1(u)\chi_{_{\mathcal
U}}(u)=\rho_2(u)\chi_{_{\mathcal U}}(u)\ a.e.$. We add and subtract
$\sum_i\int_{\mathcal U_i}h_{m_2}\left(\rho_2^{\mathcal
U_i}(u)\right)du$ and finally we obtain
\begin{equation}
S_2(\rho_1,\rho_2)-S_1(\rho_1)=\int_{\mathcal
U^c}h_{m_2}(\rho_2(u))du+ \sum_i\int_{\mathcal
U_i}h_{m_2}\left(\rho_2^{\mathcal U_i}(u)\right)\geq 0\ .
\label{gengivkan}
\end{equation}
We used the fact that given an interval $\mathcal V$ and a density
$\rho$
$$
\int_{\mathcal
V}\left(h_{m_2}(\rho(u))-h_{m_1}(\rho(u))\right)du=\left(\int_{\mathcal
V}\rho(u)
du\right)\log\frac{m_1(1-m_2)}{m_2(1-m_1)}+\log\frac{1-m_1}{1-m_2}\
,
$$
depends only on $m_1$, $m_2$ and the total mass $\int_{\mathcal
V}\rho(u) du$. This imply that
$$
\int_{\mathcal
U_i}\left(h_{m_2}(\rho_1(u))-h_{m_1}(\rho_1(u))\right)du=
\int_{\mathcal U_i}\left(h_{m_2}(\rho_1^{\mathcal
U_i}(u))-h_{m_1}(\rho_1^{\mathcal U_i}(u))\right)du \ ,
$$
due to the fact that $\int_{\mathcal U_i}\rho_1(u) du=\int_{\mathcal
U_i}\rho_1^{\mathcal U_i}(u) du$. The right hand side of
\eqref{gengivkan} can be zero if and only if
$\rho_2(u)\chi_{_{\mathcal U^c}}(u)=m_2\chi_{_{\mathcal U^c}}(u)\
a.e.$ and $\rho_2^{\mathcal U_i}(u)\chi_{_{\mathcal U_i}}(u)=m_2
\chi_{_{\mathcal U_i}}(u)\ a.e.$ for any $i$. This happens if and
only if $\rho_2$ is constructed as in \eqref{ggg}. Let us show this
fact.

Given two intervals $\widetilde{\mathcal V}_i$ and
$\widetilde{\mathcal V}_j$ then they are either disjoint or one
contained inside the other. This follows from the following
statement: if $w_i^l\not \in [v_j^r,v_i^l)$ then $w_i^l\in
(v_i^r,w_j^l)$. To prove the statement observe that, by definition
of $w_i^l$, for any $u \in (w_i^l,v_i^r)$ we have
$$
\int_{w_i^l}^u\left(m_2-\rho_1(z)\right)
dz+\int^{v_i^r}_u\left(m_2-\rho_1(z)\right) dz=0
$$
and moreover the second integral is strictly negative. As a
consequence for any $u \in (w_i^l,v_i^r)$ it holds
\begin{equation}
\int_{w_i^l}^u\left(m_2-\rho_1(z)\right)dz>0\ . \label{caccodrillo}
\end{equation}
If $w_i^l\not \in [v_j^r,v_i^l)$ then
$\int_{w_i^l}^{v_j^r}\left(m_2-\rho_1(z)\right)dz>0$. Clearly we
have also $\int_{u}^{v_j^r}\left(m_2-\rho_1(z)\right)dz<0$ for any
$u\in (w_j^l,v_j^r)$. As a consequence we deduce that there exists
an $u\in (w_i^l,v_j^l)$ such that
$\int_{u}^{v_j^r}\left(m_2-\rho_1(z)\right)dz=0$ and this implies
the above statement.

We consider the subfamily of intervals $\widetilde{\mathcal
V}_{i_k}$ composed by the intervals $\widetilde{\mathcal V}_i$ that
are not contained inside intervals $\widetilde{\mathcal V}_j$ with
$j\neq i$. Note that $\cup_k \widetilde{\mathcal V}_{i_k}=\cup_i
\widetilde{\mathcal V}_i$.

The fact that \eqref{ggg} is a minimizer of \eqref{gengivkan} and
consequently also of \eqref{pierpo} follows from the fact that
inequality \eqref{caccodrillo} is equivalent to inequality
$F_{m_2}^{\widetilde{\mathcal V}_{i_k}}\geq
F_{\rho_1}^{\widetilde{\mathcal V}_{i_k}}$ with the equality sign
holding only at the boundary of $\widetilde{\mathcal V}_{i_k}$. This
implies the fact that $\widehat{F}_{\rho_1}^{\widetilde{\mathcal
V}_{i_k}}$ coincides with $F_{m_2}^{\widetilde{\mathcal V}_{i_k}}$
and consequently $\rho_1^{\widetilde{\mathcal
V}_{i_k}}(u)\chi_{_{\widetilde{\mathcal V}_{i_k}}}(u)=m_2
\chi_{_{\widetilde{\mathcal V}_{i_k}}}(u)\ a.e.$. The fact that
\eqref{ggg} is the unique minimizer can be shown from the fact that
\eqref{gengivkan} is strictly positive for different density
profiles.
\end{proof}

The density profile $\rho_2^*$ in \eqref{ggg} is the typical total
density profile when the system is conditioned to have a density
profile $\rho_1$ of first class particles.

\vskip0.5cm

Following the steps of all the proofs presented for the TASEP you
see that the only properties of the model that we used are: the
strict convexity in $x$ of $h_m(x)$ and the the fact that
$h_m(x)\geq 0$ with equality if and only if $x=m$. Both properties
hold also for $k_m(x)$. As a consequence all the above statements
hold also for the the HAD process. Starting from theorem
\ref{bernoulli2} and proceeding as before we obtain the following
result.

\begin{Th}
\label{Th2} Let $0<m_1<m_2$ positive real numbers. Consider the
$2$-class HAD process on $\Lambda$ having respectively $[Nm_1]$
first class particles and $[Nm_2]$ total particles. When the pair
$(\underline{x}^{(1)},\underline{x}^{(2)})$ is distributed according
to the invariant measure of the process, we have that
$\pi_N^2(\underline{x}^{(1)},\underline{x}^{(2)})$ satisfies a large
deviation principle with parameter $N$ and with good rate function
$S_2(\rho_1,\rho_2)$ defined as follows. It takes the value
$+\infty$ when $(\rho_1,\rho_2)\not \in
I^{2,\uparrow}\cap\left(\mathcal M^{0}_{m_1}\times \mathcal
M^{0}_{m_2}\right)$. When $(\rho_1,\rho_2) \in
I^{2,\uparrow}\cap\left(\mathcal M^{0}_{m_1}\times \mathcal
M^{0}_{m_2}\right)$ it takes the value
\begin{equation}
S_2(\rho_1,\rho_2)= \int_{\mathcal U^c}k_{m_1}(\rho_1(u))du
+\sum_i\int_{\mathcal U_i}k_{m_1}\left(\rho_1^{\mathcal
U_i}(u)\right)+\int_{\Lambda}k_{m_2}(\rho_2(u))du\ ,
\end{equation}
where $ \mathcal U_i$ are disjoint closed intervals and $\mathcal
U=\cup_i \mathcal U_i$ is a.e. equivalent to the set
$$
\left\{u\in \Lambda : \ \rho_1(u)=\rho_2(u) \right\}\ .
$$
The symbol $\rho_1^{\mathcal U_i}$ is defined before the statement
of lemma \ref{gigolino}.
\end{Th}
The rate functional $S_2$ is non negative and is zero if and only if
$\psi_1(u)=m_1 \ a.e.$ and $\psi_2(u)=m_2 \ a.e.$ and corresponding
$\rho_2(u)=m_2 \ a.e.$ and $\rho_1(u)=C_{m_2}[m_1](u)=m_1 \ a.e.$.

Also in this case the rate functional $S_2$ is not convex.

The typical density $\rho_1^*$ of first class particles when you
condition the system to have a total density $\rho_2$ is given by
\eqref{fff}.

The typical total density profile $\rho_2^*$ for the system
conditioned to have a density profile $\rho_1$ of first class
particles is given by \eqref{ggg}.

\section{LDP for multiclass processes}
\label{ldpk}

Theorem \ref{bernoulli} is also immediately generalized to the case
of product of $k$ measures. Consider $0<m_1\leq m_2\leq\leq \dots
\leq m_k\leq 1$ real numbers. The family of measures
\begin{equation}
\left(\nu_N^{[Nm_1]}\times \cdots \times\nu_N^{[Nm_k]}\right)\circ
\left(\pi_N^{k}\right)^{-1}\nonumber
\end{equation}
satisfies a LDP with parameter $N$ and with good rate functional
given by
\begin{equation}
\left\{
\begin{array}{cc}
\sum_{i=1}^k\int_{\Lambda} du \ h_{m_i}(\rho_i(u)) & if \ \
\rho_i\in \mathcal M^{0,b}_{m_i}\ ,\\
+\infty & otherwise\ .\\
\end{array}
\right. \nonumber
\end{equation}
We are interested in proving a LDP for the empirical measures of the
invariant measures of the k-class TASEP. This means that we are
interested in proving a LDP for the sequence of measures
\begin{equation}
\left[\left(\nu_N^{[Nm_1]}\times \dots \times
\nu_N^{[Nm_k]}\right)\circ \left(\mathbb
C_{k}\right)^{-1}\right]\circ \left(\pi_N^k\right)^{-1}\ , \nonumber
\end{equation}
that due to identity \eqref{supercommuta} coincides with the
sequence of measures
\begin{equation}
\left[\left(\nu_N^{[Nm_1]}\times \dots \times
\nu_N^{[Nm_k]}\right)\circ \left(\pi_N^{k}\right)^{-1}\right]\circ
\left(\mathbb C_k\right)^{-1}\ . \nonumber
\end{equation}
As in the previous section we can apply the contraction principle
obtaining the following result.

\begin{Th}
Let $0<m_1\leq \dots \leq m_k\leq 1$ be real numbers. Consider the
$k$-class TASEP on $\mathbb Z_N$ with $[Nm_i]$ particles of class
$\leq i$. When $(\eta_1,\dots ,\eta_k)$ is distributed according to
the invariant measure of the process, we have that
$\pi_N^k(\eta_1,\dots ,\eta_k)$ satisfies a LDP with parameter $N$
and good rate function $S_k(\rho_1, \dots ,\rho_k)$ given by
\begin{equation}
S_k(\rho_1,\dots ,\rho_k)=\inf_{\left\{(\psi_1,\dots ,\psi_k): \
\psi_i\in \mathcal M_{m_i}^{0,b},\ \mathbb C_k(\psi_1, \dots
,\psi_k)=(\rho_1, \dots
,\rho_k)\right\}}\left(\sum_{i=1}^k\int_{\Lambda}h_{m_i}(\psi_i(u))du\right)\
,\label{ultimo}
\end{equation}
with the convention that the infimum over an empty set is $+\infty$.
\end{Th}

Remember that in \eqref{ultimo} we are indicating with $\psi_i$ both
the measure and the corresponding density. The functional $S_k$ is
nonnegative and zero if and only if $\rho_i(u)=m_i \ a.e.$, moreover
it takes the value $+\infty$ on
$\left[I^{k,\uparrow}\cap\left(\mathcal M_{m_1}^{0,b}\times \dots
\times \mathcal M_{m_k}^{0,b}\right)\right]^c$.

We obtained an interesting variational problem that will not be
studied in this paper. We remark only the following fact. Existence
of a minimizer for \eqref{ultimo} can be proved using the same
strategy as in theorem \ref{Th1}. Uniqueness of the minimizer is not
guaranteed as in the case of $2$-class models. We have in fact that
when $k>2$ the set
\begin{equation}
\left\{(\psi_1,\dots ,\psi_k): \ \psi_i\in \mathcal M_{m_i}^{0,b},\
\mathbb C_k(\psi_1,\dots \psi_k)=(\rho_1,\dots ,\rho_k)\right\}
\label{inatteso}
\end{equation}
is not necessarily a convex set. This follows from the following
example. Let $\epsilon$ be any positive real number $<\frac{1}{8}$
and consider the measures defined from the following densities
\begin{equation}
\left\{
\begin{array}{l}
\psi_1=2\chi_{[\frac{1}{8},\frac{1}{8}+\epsilon]}\ ,\\
\psi_2=4\chi_{[0,\epsilon]}+4\chi_{[\frac{7}{8},\frac{7}{8}+\epsilon]}\ ,\\
\psi_3=4\chi_{[\frac{1}{4},\frac{1}{4}+\epsilon]}+8\chi_{[\frac{1}{2},\frac{1}{2}+\epsilon]}\
,
\end{array}
\right.\nonumber
\end{equation}
and
\begin{equation}
\left\{
\begin{array}{l}
\widetilde{\psi}_1=2\chi_{[\frac{5}{8},\frac{5}{8}+\epsilon]}\ ,\\
\widetilde{\psi}_2=4\chi_{[\frac{3}{8},\frac{3}{8}+\epsilon]}+4\chi_{[\frac{3}{4},\frac{3}{4}+\epsilon]}\ ,\\
\widetilde{\psi}_3=4\chi_{[\frac{1}{4},\frac{1}{4}+\epsilon]}+8\chi_{[\frac{1}{2},\frac{1}{2}+\epsilon]}\
.
\end{array}
\right.\nonumber
\end{equation}
We have that
$$\mathbb C_3(\psi_1,\psi_2,\psi_3)=\mathbb
C_3(\widetilde{\psi}_1,\widetilde{\psi}_2,\widetilde{\psi}_3)=(\rho_1,\rho_2,\rho_3)\
,
$$
where the measures $\rho_i$ have densities a.e. equal to
\begin{equation}
\left\{
\begin{array}{l}
\rho_1=4\chi_{[\frac{1}{4},\frac{1}{4}+\frac{\epsilon}{2}]}\ ,\\
\rho_2=4\chi_{[\frac{1}{4},\frac{1}{4}+\epsilon]}+
8\chi_{[\frac{1}{2},\frac{1}{2}+\frac{\epsilon}{2}]}\ ,\\
\rho_3=4\chi_{[\frac{1}{4},\frac{1}{4}+\epsilon]}+8\chi_{[\frac{1}{2},\frac{1}{2}+\epsilon]}\
.
\end{array}
\right.\nonumber
\end{equation}
It is easy to check that for the convex combination
$$
(\phi_1,\phi_2,\phi_3)=\frac{1}{2}(\psi_1,\psi_2,\psi_3)+\frac{1}{2}
(\widetilde{\psi}_1,\widetilde{\psi}_2,\widetilde{\psi}_3)
$$
we have $\mathbb
C_3(\phi_1,\phi_2,\phi_3)\neq(\rho_1,\rho_2,\rho_3)$.

\vskip0.5cm

The same kind of result is easily derived also for the HAD process
starting from theorem \ref{bernoulli2}.

\begin{Th}
Let $0<m_1\leq \dots \leq m_l$ be real numbers. Consider the
$l$-class HAD process on $\Lambda$ with $[Nm_i]$ points of class
$\leq i$. When $(\underline{x}^{(1)},\dots ,\underline{x}^{(l)})$ is
distributed according to the invariant measure of the process, we
have that $\pi_N^l(\underline{x}^{(1)},\dots ,\underline{x}^{(l)})$
satisfies a LDP with parameter $N$ and with good rate function
$S_l(\rho_1, \dots ,\rho_l)$ given by
\begin{equation}
S_l(\rho_1,\dots ,\rho_l)=\inf_{\left\{(\psi_1,\dots ,\psi_l): \
\psi_i\in \mathcal M_{m_i}^{0}, \  \mathbb C_l(\psi_1, \dots
,\psi_l)=(\rho_1, \dots
,\rho_l)\right\}}\left(\sum_{i=1}^l\int_{\Lambda}k_{m_i}(\psi_i(u))du\right)
\end{equation}
with the convention that the infimum over an empty set is $+\infty$.
\end{Th}
The functional $S_l$ is nonnegative and zero if and only if
$\rho_i(u)=m_i \ a.e.$. Moreover it takes the value $+\infty$ on
$\left[I^{l,\uparrow}\cap\left(\mathcal M_{m_1}^{0}\times \dots
\times \mathcal M_{m_l}^{0}\right)\right]^c$.

\vskip0.5cm

Still using the contraction principle we obtain the following
identity, valid both for the TASEP and the HAD process, whose study
from the variational point of view seems to be interesting
\begin{equation}
\inf_{\left\{\rho_j\right\}}S_k(\rho_1,\dots
,\rho_k)=S_{k-1}(\rho_1, \dots ,\widehat{\rho}_j,\dots, \rho_k)\ .
\label{contrall}
\end{equation}
With the symbol $\widehat{\rho}_j$ we indicate the fact that the
measure $\rho_j$ is missing. From the microscopic point of view
identity \eqref{contrall} derives from the fact that if you change
the class of $j$-class particles to $j+1$ the dynamics that you
observe is the one of a $(k-1)$-class process.

We derive now a recursive relation. We derive it for the TASEP but
it holds also for the HAD process. To simplify notations the fact
that all the measures $\rho_i$, $\psi_i$ and $\phi_i$ involved are
absolutely continuous, have bounded densities and have fixed total
mass will be understood. Also the fact that $(\rho_1,\dots
,\rho_k)\in I^{k,\uparrow}$ is understood. We can write
\eqref{ultimo} as
\begin{eqnarray}
& S_k(\rho_1,\dots ,\rho_k)=\int_{\Lambda}h_{m_k}(\rho_k(u))du
\nonumber \\
& +\inf_{\left\{(\psi_1,\dots ,\psi_{k-1}): \ \mathbb C_k(\psi_1,
\dots ,\psi_{k-1},\rho_k)=(\rho_1, \dots
,\rho_k)\right\}}\left(\sum_{i=1}^{k-1}\int_{\Lambda}h_{m_i}(\psi_i(u))du\right)\
, \label{veroultimo}
\end{eqnarray}
because if $C_k(\psi_1, \dots ,\psi_{k})=(\rho_1, \dots ,\rho_k)$
then $\psi_k=\rho_k$. Let us call $(\phi_1,\dots
,\phi_{k-1}):=\mathbb C_{k-1}(\psi_1,\dots ,\psi_{k-1})$. We have
then that $(\phi_1,\dots ,\phi_{k-1})\in A_{\rho_k}^{k-1}$, where
$$
A^{k-1}_{\rho_k}:=\left\{(\phi_1, \dots ,\phi_{k-1})\in
I^{k-1,\uparrow}: \ C_{\rho_k}[\phi_i]=\rho_i, \ i=1,\dots
,k-1\right\}\ .
$$
Then equation \eqref{veroultimo} becomes
\begin{eqnarray*}
& S_k(\rho_1,\dots ,\rho_k)=\int_{\Lambda}h_{m_k}(\rho_k(u))du
\nonumber \\
& +\inf_{\left\{(\phi_1,\dots ,\phi_{k-1})\in
A_{\rho_k}^{k-1}\right\}}\inf_{\left\{(\psi_1,\dots ,\psi_{k-1}): \
\mathbb C_{k-1}(\psi_1, \dots ,\psi_{k-1})=(\phi_1, \dots
,\phi_{k-1})\right\}}\left(\sum_{i=1}^{k-1}\int_{\Lambda}h_{m_i}(\psi_i(u))du\right)
\end{eqnarray*}
that finally becomes the following recursive relation
$$
S_k(\rho_1,\dots
,\rho_k)=\int_{\Lambda}h_{m_k}(\rho_k(u))du+\inf_{\left\{(\phi_1,\dots
,\phi_{k-1})\in A_{\rho_k}^{k-1}\right\}}S_{k-1}(\phi_1,\dots
,\phi_{k-1})\ .
$$

\section{Acknoledgments}
I thank L. Bertini and A. Faggionato for useful discussions and
suggestions.

 \end{document}